\input amstex
\documentstyle{amsppt}
\magnification=1200
%\hcorrection{0.7in}
%\vcorrection{-0.7in}
\NoRunningHeads
\NoBlackBoxes
\topmatter
\title Collapsed $5$-manifolds with pinched positive sectional curvature
\endtitle
\author Fuquan Fang \footnote{Supported by NSFC Grant
19925104 of China, 973 project of Foundation Science of China
\hfill{$\,$}}\,\& \,Xiaochun Rong\footnote{Supported partially by
NSF Grant DMS 0203164 and a research found from Capital normal
university \hfill{$\,$}}
\endauthor
\address Mathematics Department, Capital Normal University,
Beijing, P.R.C.
\newline
.\hskip3mm Nankai Institute of Mathematics, Nankai University,
Tianjin 300071, P.R.C.
\endaddress
\email ffang\@nankai.edu.cn
\endemail
\address Mathematics Department, Rutgers University, New Brunswick,
NJ 08903, U.S.A
\newline
.\hskip3mm Mathematics Department, Capital Normal University,
Beijing, P.R.C.
\endaddress
\email rong\@math.rutgers.edu
\endemail
\abstract Let $M$ be a closed $5$-manifold of pinched curvature
$0<\delta\le \text{sec}_M\le 1$. We prove that $M$ is homeomorphic
to a spherical space form if $M$ satisfies one of the following
conditions: (i) $\delta =1/4$ and the fundamental group is a
non-cyclic group of order $\ge C$, a constant. (ii) The center of
the fundamental group has index $\ge w(\delta)$, a constant
depending on $\delta$. (iii) The ratio of the volume and the
maximal injectivity radius is $<\epsilon(\delta)$. (iv) The volume
is less than $\epsilon(\delta)$ and the fundamental group
$\pi_1(M)$ has a center of index at least $w$, a universal
constant, and $\pi_1(M)$ is either isomorphic to a spherical
$5$-space group or has an odd order.
\endabstract
\endtopmatter
\document

\head 0. Introduction
\endhead

\vskip4mm

The sphere theorem asserts that if a manifold $M$ admits a quarter
pinched metric, $\frac 14<\text{sec}_M\le 1$, then its universal
covering space is homeomorphic to a sphere. A natural problem is
to determine whether $M$ is homeomorphic to a spherical space
form, $S^n/\Gamma$, $\Gamma\subset O(n+1)$? We will call $\Gamma$
a {\it spherical $n$-space group}. This problem is (wildly) open
in all odd-dimension $n\ge 5$ (cf. [Ha], [GKR], [IHR]). Clearly, a
positive answer implies that the fundamental group $\pi_1(M)$ is
isomorphic to a spherical space group and the $\pi_1(M)$-action on
the universal covering is conjugate to a linear action. A subtlety
is that neither of these holds without a positive curvature
condition: in every odd dimension $n\ge 5$ ([Ha]),

\noindent (0.1) There are infinitely many non-spherical space
groups acting freely on an $n$-sphere.

\noindent (0.2) There are infinitely many distinct free actions on
an $n$-sphere by a spherical space group which do not conjugate to
any linear action.

\noindent Obviously, the quotient manifolds in (0.1) and (0.2) are
not homeomorphic to any spherical space form.

In this paper, as a first step we investigate the case of
dimension $5$. We give positive answers (Theorems A-D) for a
$\delta$-pinched $5$-manifold $M$ whose fundamental group is
not small (equivalently, whose volume is small).
In particular, we
rule out (0.1) and (0.2) in our circumstances via studying certain
symmetry structure on $M$ discovered by Cheeger-Fukaya-Gromov
([CFG], [Ro1]).

We now begin to state the main results in this paper.

\vskip2mm

\proclaim{Theorem A}

Let $M$ be a closed $5$-manifold with $\frac 14 <\text{sec}_M\le
1$. If the fundamental group $\pi_1(M)$ is a non-cyclic group of
order $\ge C$ (a constant), then $M$ is homeomorphic to a
spherical space form.
\endproclaim

\vskip2mm

For any $\delta$-pinching, we can generalize Theorem A (whose assumption is slightly stronger than the case of $\delta=1/4$ in Theorem B).

\vskip2mm

\proclaim{Theorem B}

Let $M$ be a closed $5$-manifold with $0<\delta\le \text{sec}_M\le
1$. If the center of $\pi_1(M)$ has index $\ge w(\delta)$ (a
constant depending on $\delta$), then $M$ is homeomorphic to a
spherical space form.
\endproclaim

\vskip2mm

Note that $M$ satisfies that $\text{diam}(M)\le \pi/\sqrt \delta$
(Bonnet theorem) and $\text{vol}(M)\le
\text{vol}(S^5_\delta)/w(\delta)<<1$ (volume comparison), where
$S^5_\delta$ denotes the sphere of constant curvature $\delta$.
According to [CFG] (cf. [Ro1]), $M$ admits local compatible
isometric $T^k$-actions with $k\ge 1$ of some nearby positively curved metric
(details will be given shortly).  In our circumstance, we show
that $k\ge 2$, and thus the universal covering of $M$ is
diffeomorphic to a sphere ([Ro2]). The main work in the proofs
of Theorems A and B is to show, using
the local symmetry structure, that $\pi_1(M)$ is isomorphic to a
spherical space group and the $\pi_1(M)$-action is conjugate to a
linear one (see Theorem E, compare to (0.1) and (0.2)).

Observe that any $\delta$-pinched $5$-manifold satisfies that
$\frac{\text{vol}(M)}{\max \{\text{injrad}(M,z)\}}\le d(\delta)$ (by the
Cheeger's lemma, [CE]). When the ratio is small, we find that
the above isometric $T^k$-action satisfies $k\ge 2$ (note that
because the Riemannian universal covering is a sphere ([Ro2],
$|\pi_1(M)|$ is propositional to $\text{vol}(M)$, [PRT]).

\vskip2mm

\proclaim{Theorem C}

For $0<\delta\le 1$, there exists a small number,
$\epsilon(\delta)>0$, such that if a closed $5$-manifold $M$ satisfies
$$ 0<\delta\le \text{sec}_M\le 1,\qquad \qquad \frac{\text{vol}
(M)}{\max\text{injrad }(M,x)}<\epsilon(\delta),$$ then $M$ is
homeomorphic to a spherical space form.
\endproclaim

\vskip2mm

Note that Theorems B and C do not hold if one replaces the
condition ``$\delta>0$'' with ``$\delta\ge 0$'' (see Example 2.7).
We intend to discuss the classification with ``$\delta\ge 0$''
elsewhere.

Consider a collapsed $\delta$-pinched $5$-manifold $M$ close in
the Gromov-Hausdorff distance to a metric space $X$ of dimension
$4$ (equivalently, $k=1$). A new trouble is to determine the
topology of the universal covering space, or equivalently the
topology of $X$. This looks quite difficult; it seems to require a
classification of positively curved $4$-manifolds in the case when
$X$ is smooth (and thus the fundamental group of $M$ is cyclic).
The following result says that $M$ is homeomorphic to a spherical
space form when $\pi_1(M)$ is certain non-cyclic group.

\vskip2mm

\proclaim{Theorem D}

For $0<\delta\le 1$, there exists $\epsilon(\delta)>0$ such
that if a closed $5$-manifold $M$ satisfies
$$ 0<\delta\le \text{sec}\le 1,\qquad \qquad \text{vol}(M)
<\epsilon(\delta),$$ then $M$ is homeomorphic to a spherical space
form, provided $\pi_1(M)$ has a center of index at least $w>0$, a
constant (independent of $\delta$), and $\pi _1(M)$ is a spherical
$5$-space group or $|\pi _1(M)|$ is odd.
\endproclaim

\vskip2mm

We mention that there are infinitely many spherical $5$-space
forms satisfying the conditions of Theorems A-D (see Example 2.6).

A natural question is when $M$ in Theorems A-D is diffeomorphic to
a spherical $5$-space form? We mention the following: a spherical
$5$-space form $S^5/\Gamma$ admits exactly one or two different
smooth structures depending on $|\Gamma|$ odd or even ([KS],
[Wa]). Moreover, both smooth structures may allow a non-negatively
curved metric (e.g., there are exactly four smooth manifolds
homotopy equivalent to $\Bbb RP^5$, all of them admit metrics of
non-negative sectional curvature ([GZ]), and two of them are not
homeomorphic to each other).

A question of Yau ([Yau]) is if in a given homotopy type contains
at most finitely many diffeomorphism types that can
support a metric of positive sectional curvature. A positive
answer is known only in dimensions $2$ and $3$ ([Ha]). When
restricting to the class of pinched metrics, positive answers are
known in even dimensions and the class of manifolds
(odd-dimensions) with finite second homotopy groups ([FR1], [PT]).

In view of the above, Theorem B has the following corollary.

\vskip2mm

\proclaim{Corollary 0.3}

Let $M$ be a close $\delta$-pinched $5$-manifold. Then the
homotopy type of $M$ contains at most $c(\delta)$ many
diffeomorphism types that support a $\delta$-pinched metric,
provided $\pi_1(M)$  has a center with index $\ge w(\delta)$.
\endproclaim

\vskip2mm

As mentioned earlier, our approach to Theorems A-D is based on the
fibration theorem of Cheeger-Fukaya-Gromov on collapsed manifolds
with bounded sectional curvature and diameter ([CFG], [CG1,2]). In
our circumstances, the fibration theorem asserts that there is a
constant $v(n,\delta)>0$ such that if a $\delta$-pinched
$n$-manifold $M$ has a volume less than $v(n,\delta)$, then $M$
admits a pure F-structure all whose orbits are of positive
dimensions. By the Ricci flows technique, one can show that there
is invariant metric which is at least $\delta/2$-pinched ([Ro1]).

In the case of a finite fundamental group, the notion of a pure
F-structure is equivalent to that of a {\it $\pi_1$-invariant
torus $T^k$-action} on a manifold $M$, which is defined by an
effective $T^k$-action on the universal covering space $\tilde M$
of $M$ such that it extends to a $T^k\rtimes
_\rho\pi_1(M)$-action, where $\rho: \pi_1(M)\to \text{Aut}(T^k)$
is a homomorphism from the fundamental group to the automorphism
group of $T^k$. Clearly, the $T^k$-action on $\tilde M$ is the
lifting of a $T^k$-action on $M$ if and only if $\rho$ is trivial
or equivalently, the $T^k$-action and the $\pi_1(M)$-action
commute. Hence, the notion of a $\pi_1$-invariant torus action
generalizes that of a global torus action and the $T^k$-orbit
structure on $\tilde M$ projects onto $M$ so that each orbit is a
flat submanifold.

Consider $M$ as in Theorems A-D. In view of the above, we may
assume that $M$ admits a $\pi_1$-invariant isometric $T^k$-action
($k\ge 1$).

\vskip2mm

\proclaim{Theorem E}

Let $M$ be a closed $5$-manifold of positive sectional curvature.
If $M$ admits a $\pi_1$-invariant isometric $T^k$-action with $k>1$,
then $M$ is homeomorphic to a spherical space form.
\endproclaim

\vskip2mm

Theorem E is known in the following special cases: $M$ (itself)
admits an isometric $T^3$-action ([GS]) or $M$ is simply connected
([Ro2]).

We show that under the assumptions of Theorems B and C, $k>1$
and thus Theorem E implies Theorems B and C. In the case $k=1$, the
$T^1$-action on the universal covering $\tilde M$ of $M$ is free
if $\tilde M$ is a sphere and if $\pi_1(M)$ is not cyclic. We then
complete the proof of Theorem A by proving the following
topological result: Let a finite group $\Gamma$ act freely on
$S^5$. If $S^5$ admits a free $\Gamma$-invariant $T^1$-action such
that the induce $\Gamma$-action on $S^5/T^1$ is pseudo-free, then
$S^5/\Gamma$ is homeomorphic to a spherical space form (see
Proposition 1.4).

In view of the above, Theorem D follows from

\vskip2mm

\proclaim{Theorem F}

Let $M$ be a closed $5$-manifold of positive sectional curvature
which admits a $\pi_1$-invariant fixed point free isometric
$T^1$-action. Then the universal covering $\tilde M$ is
diffeomorphic to $S^5$, provided $\pi_1(M)$ has a center of index
$\ge w$, and $\pi _1(M)$ is a spherical
$5$-space group or $|\pi _1(M)|$ is odd.
\endproclaim

\vskip2mm

It is worth to point it out that {\it every} spherical $5$-space
form admits a $\pi_1$-invariant isometric $T^3$-action and a free
isometric $T^1$-action and a $\pi_1(M)$-invariant
isometric $T^2$-action ([Wo], p.225).

We would like to put Theorems E and F in a little perspective. In
the study of positive sectional curvature, due to the obvious
ambiguity the class of positively curved manifolds with (large)
symmetry has frequently been a focus of the investigations.
According to K. Grove, this also serves as a strategy of searching
for new examples and obstructions. There has been significant
progress in the last decade on classification of simply connected
manifolds with large symmetry rank (the rank of the isometry
group), cf. [GS], [FR1,2], [HK], [Ro1-3], [Wi1,2]. However, not
much is known for non-simply connected manifolds with large
symmetry rank. Theorems E and F may be treated as an attempt in
this direction.

We now give an outline of the proof of Theorems E and F.

By the compact transformation group theory ([Bre]), the topology
of a $T^k$-space $M$ is closely related to that of the singular
set (the union of all non-principal orbits) and the orbit space
$M/T^k$. In the presence of an invariant metric of positive
sectional curvature, the singular set and the orbit space are very
restricted, and this is the ultimate reason for a possible
classification. In our proofs, we will thoroughly investigate the
structure of the singularity and the orbit space.

The proof of Theorem E divides into two situations: $k=3$ (Theorem
3.1) and $k=2$, and the main work is in the case of $k=2$. When
$k=2$, we first prove Theorem E at the level of fundamental groups
(Theorem 4.1). Then we divide the proof into two cases: the
$\pi_1$-invariant $T^2$-action is pseudo-free (section 5) and not
pseudo-free (section 6). In the former case, we study the
$T^2$-action on the universal covering space, $\tilde M\simeq
S^5$, which has a singular set, $\Cal S$, consisting of three
isolated circle orbits. It suffices to show that the
$(\pi_1(M),T^2)$-bundle, $\tilde M-\Cal S\to (\tilde M-\Cal
S)/T^2$, is conjugate to a standard linear model from spherical
space forms. This can be done following [FR2] if one can assume
that the orbit space, $\tilde M^*=\tilde M/T^2$, is a homeomorphic sphere.
The problem is that we only know that $X$ is a homotopy
$3$-sphere. We overcome this difficulty by combining the above
with tools from the $s$-cobordism theory in dimension $5$. In the
non-pseudo-free case, the singular set has dimension $3$, and by
analyzing the singular structure, we are able to view $M$ as a
gluing of standard pieces.

In the proof of Theorem F, by studying the induced
$\pi_1(M)$-action on $\tilde M^*=\tilde M/T^1$ (which is not trivial because
$\pi_1(M)$ is not cyclic), we bound from above the Euler
characteristic of $\tilde M^*$ via the technique of $q$-extent
(Proposition 9.5, cf. [Gr], [Ya]). With this constraint, we show
that the condition on the fundamental groups implies that $\tilde
M$ is a homology sphere, the $T^1$-action is free and the
$\pi_1(M)$-action on $\tilde M^*$ is pseudo-free (note that the
standard free linear $T^1$-action on $S^5$ is preserved by any
spherical $5$-space group, [Wo]). By employing results in [HL] and
[Wi1,2] on pseudo-free actions by finite groups on a homeomorphic
complex projective plane, we show that the $\pi_1(M)$-action is
homeomorphically conjugate to a linear action.

\vskip2mm

\remark{Remark \rm 0.5}

By [AM], the $\frac 14$-pinching in Theorem A may be
replaced by a slightly weaker pinching constant
$\frac 14-\epsilon$.
\endremark

\vskip2mm

\remark{Remark \rm 0.6}

Theorem F actually holds with a weak restrictions on $\pi _1(M)$ (see Section 3).
\endremark

\vskip2mm

The rest of the paper is organized as follows:

Part I. Proofs of Theorems A-D by assuming Theorems E and F.

\vskip 2mm

$\S1$. Proof of Theorem A by assuming Theorem E.

$\S2$. Proof of Theorems B and C by assuming Theorem E.

$\S3$. Proof of Theorem D by assuming Theorems E and F.

\vskip 2mm

Part II. Proofs of Theorems E and F

$\S4$. Preparations

$\S5$. Proof of Theorem E for $k=3$.

$\S6$. Proof of Theorem E at the level of fundamental groups.

$\S7$. Proof of Theorem E for pseudofree $T^2$-actions.

$\S8$. Completion of the proof of Theorem E.

$\S9$. Proof of  Theorem F.

\vskip 2mm

Appendix: The $s$-cobordism theory.

\vskip5mm

\noindent {\bf Acknowledgment:} The first author would like to thank Michael
Mccooey and Ian Hambleton for invaluable discussions concerning group actions.

\vskip10mm

\head {\bf Part I. Proof of Theorems A-D by assuming Theorems E and F}\endhead

\vskip4mm

\head 1. Proof of Theorem A by Assuming Theorem E
\endhead

\vskip4mm

Let $M$ be as in Theorem A. Then the universal covering space
$\tilde M$ is homeomorphic to $S^5$ (the sphere theorem). Because
the volume of $M$ is small, $\tilde M$ admits a $\pi_1$-invariant
isometric $T^k$-action (Theorems 1.1 and 1.2). We will use this
structure to show that $M$ is homeomorphic to a spherical space
form (compare to (0.1) and (0.2)). By Theorem E, we may assume
that $k=1$. Because $\pi_1(M)$ is non-cyclic, we show that the
isometric $T^1$-action must be free and commute with the
$\pi_1(M)$-action such that the induced $\pi_1(M)$-action on
$\tilde M^*$ is pseudo-free (Lemma 1.10). These properties are
all that is required to prove, based on a result (Theorem 1.5,
c.f. [HL], [Wil1,2]), that $\pi_1(M)$ is isomorphic to a spherical
$5$-space group (Lemma 1.7) and the $\pi_1(M)$-action is conjugate
to a linear action (Proposition 1.4).

\vskip4mm

\subhead a. $\pi_1$-invariant isometric $T^k$-actions
on collapsed manifolds
\endsubhead

\vskip4mm

According to [CFG], if a closed $n$-manifold $M$ with bounded
curvature and diameter has a small volume, then $M$ admits a pure
nilpotent Killing structure whose orbits are infra-nilmanifolds.
When the fundamental group of $M$ is finite, an orbit is the
actually a flat manifold, and the nilpotent Killing structure is
equivalent to a $\pi_1$-invariant almost isometric $T^k$-action
([Ro1]). The $\pi_1$-invariance  implies that the $T^k$-orbits on
the universal covering descend to $M$, also denoted by $T^k(x),
x\in M$. A $T^k$-orbit on $M$ is called regular, if it has a
tubular neighborhood in which $T^k$-orbits form a fiber bundle.
Let $\Cal S$ denote the set of all non-regular orbits. Then
$M-\Cal S$ is an open dense subset. For a small number $\eta>0$,
let $U_{-\eta}(\Cal S)=\{x\in M:\, d(x,\Cal S)>\eta\}$.

\vskip2mm

\proclaim{Theorem 1.1 ([CFG])}

Given $n, d>0$, there exist constants, $\epsilon(n,d), c(n)>0$,
such that if a closed $n$-manifold
$M$ of finite fundamental group satisfying
$$|\text{sec}_M|\le 1,\qquad \text{diam }(M)\le d,\qquad \text{vol }
(M)<\epsilon\le \epsilon(n,d),$$
then $M$ admits a $\pi_1$-invariant $T^k$-action satisfying

\noindent (1.1.1) Every $T^k$-orbit has a positive dimension and diameter $<\epsilon$.

\noindent (1.1.2) Any orbit in $U_{-\eta}(\Cal S)$ has a second
fundamental form, $|II|\le c(n)\eta^{-1}$.

\noindent (1.1.3) There is a $T^k$-invariant metric of sectional
curvature bounded by one which is $\epsilon$-close to the
original metric in $C^1$-norm.
\endproclaim

\vskip2mm

Note that one can always assume a small constant $\eta=\eta(n)>0$
such that $U_{-\eta}\ne \emptyset$ (if $U_{-\eta}=\emptyset$, then
$\text{diam}(M)<\eta$ and thus $M$ is almost flat ([Gr1]). Then
$\Cal S=\emptyset$ and therefore $U_{-\eta}=M$, a contradiction).
Using the Ricci flows technique ([Ha]), one can obtain a nearby metric
as follows:

\vskip2mm

\proclaim{Theorem 1.2 ([Ro1])}

Let $(M,g)$ satisfy the conditions of Theorem 1.1. For $\epsilon>0$,
there is a $T^k$-invariant metric $g_\epsilon$ such that
$$|g-g_\epsilon|_{C^1}<\epsilon,\qquad
\min \{\text{sec}_g\}-\epsilon\le \text{sec}_{g_\epsilon}\le \max
\{\text{sec}_g\}+\epsilon.$$
\endproclaim

\vskip4mm

\subhead b. Spherical $5$-space forms
\endsubhead

\vskip4mm

Before presenting proofs of Theorems A-D, it may help to
review some basic facts about a spherical $5$-space form.

Given a spherical $5$-space form, $S^5(1)/\Gamma, \Gamma\subset O(6)$,
$\Gamma$ is either cyclic or is generated by two elements,
$$A=\left(\matrix R(1/m) & 0 & 0\\
0 & R(r/m) & 0\\ 0 & 0 & R(r^2/m)\endmatrix\right),\qquad
B=\left(\matrix 0 & I & 0 \\ 0 & 0 & I\\R(3\ell/n) & 0 &
0\endmatrix\right),$$
where $R(\theta)$ denotes the standard
$2\times 2$ rotation matrix with rotation angle $2\pi \theta$,
$I$ is the $2\times 2$ identity matrix and $r, n\in \Bbb Z$ satisfy $n\equiv 0 \mod 9$, $(n(r-1),m)=1$ and $r^2+r+1\equiv 0 \mod m$ ([Wo]). Clearly,
$A^m=B^n=1$ and $BAB^{-1}=A^r$. Note that $\Gamma$ may have
linear actions $S^5(1)$ which are pairwisely non-conjugate.

Let's now construct on $S^5(1)$ a $\Gamma$-invariant isometric $T^3$-, free $T^1$- and pseudofree $T^2$-action. First, the standard
$T^3$-action on $S^5(1)$ is clearly $\Gamma$-invariant and
the free diagonal circle subgroup $T^1$-action commutes
with the $\Gamma$-action. Define a pseudo-free $T^2$-action
on $S^5(1)=\{(z_1,z_2,z_3)\in \Bbb C^3,\,\, |z_1|^2+|z_2|^2+|z_3|^2=1\}$:
$$
(e^{i\theta}, e^{i\phi})(z_1,z_2,z_3)=(e^{i(\theta+\phi)}z_1,
e^{i(\theta -2\phi)}z_2, e^{i(-2\theta+\phi)}z_3)
$$
with a principal isotropy group $\Bbb Z_3$ generated by $(e^{\frac
23 \pi i}, e^{\frac 43 \pi i})$, where $T^2=T^1_\theta\times T^1_\phi$. The holonomy representation,
$\rho: \Gamma\to \text{Aut}(T^2)$ is defined by
$$\rho(A)=\text{id}_{T^2},\qquad \rho(B)=\left(\matrix 0 & -1 \\ 1 & -1\endmatrix\right),$$
where the matrix is with respect to the standard basis $\theta,\phi$.

We conclude this section with following results (which will be
referred several times through the rest of the paper). Let $M$ be
a closed $5$-manifold of positive sectional curvature.

\vskip2mm

\noindent (1.3.1) If $M$ admits an isometric $T^3$-action, then
$M$ is diffeomorphic to a lens space ([GS]).

\noindent (1.3.2) If $M$ is simply connected and if $M$ admits an isometric $T^2$-action, then $M$ is diffeomorphic
to $S^5$ ([Ro3]).

\vskip4mm

\subhead c. A criterion of a pseudo-free linear action on
$5$-spheres
\endsubhead

\vskip4mm

From now on, we will use $S^5$ to denote a homeomorphic $5$-sphere
and $S^5(1)$ a sphere of constant curvature one. Consider a finite
group $\Gamma$ acting freely on $S^5$. As mentioned in (0.1) and
(0.2), $\Gamma$ may not be isomorphic to any spherical $5$-space
group, nor, even assuming $\Gamma$ isomorphic to a spherical
$5$-space group, the $\Gamma$-action on $S^5$ may not conjugate to
any linear action. Hence, additional
conditions are required for a $\Gamma$-action to conjugate to a
linear action.

We will give a criterion, Proposition 1.4, and use it to prove
Theorem A. Note that this criterion will be also used
in the proof of Theorem F.

Spherical space forms have been completely classified, see [Wo].
From [Wo], p225, we observe that if a finite group $\Gamma \subset
SO(6)$ acts freely on $S^5(1)$ by isometries, then $\Gamma$
commutes with a standard free linear $T^1$-action on $S^5(1)$. If,
in addition, $\Gamma$ is not cyclic, then the induced
$\Gamma$-action on $S^5/T^1$ is pseudo-free i.e., any non-trivial
element has only isolated fixed points.

The above properties are sufficient for a free $\Gamma$-action on
$S^5$ to conjugate to a linear action.

\vskip2mm

\proclaim{Proposition 1.4}

Let a finite group $\Gamma$ act freely on $S^5$.
If $\Gamma$ commutes with a free $T^1$-action on $S^5$ such
that the induced $\Gamma$-action on $S^5/T^1$ is pseudo-free, then
the $\Gamma$-action is homeomorphically conjugate to a linear action.
\endproclaim

\vskip2mm

A $G$-action is called locally linear, if each singular point has an
invariant neighborhood which is equivariantly homeomorphic to a
neighborhood of $0$ in a real representation space. In particular,
smooth actions are locally linear.

The following result (cf. [HL], [Wi1-2]) plays a crucial role in
the proof of Proposition 1.4.

\vskip2mm

\proclaim{Theorem 1.5}

Any pseudo-free locally linear
action by a finite group on a $4$-manifold homeomorphic to $\Bbb
CP^2$ is topologically conjugate to the linear action of a
subgroup of $\text{PSU} (3)$ on $\Bbb CP^2$.
\endproclaim

\vskip2mm

Note that $PSU(3)=SU(3)/\Bbb Z_3$, where $\Bbb Z_3$ is the center
of $SU(3)$. It is perhaps useful to recall some details on what
finite subgroups of $PSU(3)$ can act linearly and pseudo-freely on
$\Bbb CP^2$. By [Wi2] (also [HL]), such a group is either cyclic
$\Bbb Z_n=\langle x\rangle$, or noncyclic with a presentation
$$\{x, y: yxy^{-1}=x^r,
x^n=y^3=1, \text{ where }r^2+r+1=0\mod (n)\} \tag 1.6 $$ The
linear action of the group on $\Bbb CP^2$ is given by
$$x[z_0, z_1, z_2]=[\omega z_0, \omega ^{-r}z_1, z_2];\,\, y[z_0,
z_1, z_2]=[z_1, z_2, z_0]
$$
where $\omega=e^{\frac{2\pi i}n}$ is the $n$-th root of the unit.

Observe that the group in (1.6) is  $\Bbb Z_3\oplus \Bbb
Z_3=\langle x, y\rangle$ if $n=3$, and $n$ can not be an integral
multiple of $9$. Therefore, (1.6) can never be a $5$-dimensional
spherical space form group if $n>1$ (cf. [Wo] page 225). However,
Petrie [Pe] constructed a free action of (1.6) on $S^5$ if $n=7$
and $r=2$.

\vskip2mm

\proclaim{Lemma 1.7}

Under the assumptions of Proposition 1.4, let $\Gamma _0\subset
\Gamma$ denote the subgroup which acts trivially on $S^5/T^1$.
Then

\noindent (1.7.1) $\Gamma_0\subseteq C(\Gamma)$, the center of
$\Gamma$.

\noindent (1.7.2) $\Gamma$ is isomorphic to a subgroup of $SU(3)$.
\endproclaim

\vskip2mm

\demo{Proof} (1.7.1) Because $\Gamma_0$ acts trivially on
$S^5/T^1$, $\Gamma_0\subset T^1$, and thus $\Gamma_0\subseteq
C(\Gamma)$. Let $\Gamma_0\cong \Bbb Z_\ell$.

(1.7.2) By Theorem 1.5 and the above discussion, $\Gamma^*=\Gamma /\Gamma_0$ is either cyclic, or a noncyclic group as in (1.6). The
desired result follows in the former case, because $\Gamma$ must
be abelian and thus must be cyclic because $\Gamma$ acts freely on
$S^5$.

Let $\Gamma^*$ be a group in (1.6). Then $\Gamma^*$ has a
trivial center. Moreover, $n$ must be
coprime to $3$ because otherwise $\Gamma^*$ contains $\Bbb
Z_3\oplus \Bbb Z_3$ as a subgroup. By analyzing the central
extension $1\to \Bbb Z_\ell\to \Gamma \to \Gamma^*\to 1$
restricted on $\Bbb Z_3\oplus \Bbb Z_3$ we conclude that $\Gamma$
also contains a subgroup isomorphic to $\Bbb Z_3\oplus \Bbb Z_3$.
This is absurd since $\Gamma$ acts freely on $S^5$. Similarly, the
central extension $\Gamma$, when restricted on $\langle x\rangle $
and $\langle y\rangle $, gives rise to cyclic groups of order
$\ell n$ and $3\ell$ respectively. Hence $\Gamma$ contains $\Bbb
Z_{\ell n}$ as a cyclic subgroup of index $3$, and every Sylow
group in $\Gamma$ is cyclic.  By the Burnside Theorem (cf. [Wo]
Theorem 5.4.1, p.163), $\Gamma$ is generated by two elements $A$
and $B$ with relations
$$A^k=B^{s}=1, \text{   }BAB^{-1}=A^r$$
where $|\Gamma |=ks$, $((r-1)s, k)=1$ and $r^s\equiv 1(\text{mod
}k)$. Since $\Gamma$ contains an index $3$ normal subgroup, we
conclude that $s$ is divisible by $3$, and $\{A, B^3\}$ generates
the normal cyclic subgroup, $r^3\equiv 1(\text{mod }k)$. Note that
the center of $\langle A, B\rangle$ is generated by $B^3$ which
has order $s/3$. Hence $s=3\ell$, and $k=n$.
  Therefore, $\Gamma$ may be
realized as the subgroup of $SU(3)$ generated by the matrices
$$\left(\matrix R(1/n) & 0 & 0\\
0 & R(r/n) & 0\\ 0 & 0 & R(r^2/n)\endmatrix\right),\qquad
\left(\matrix 0 & I & 0 \\ 0 & 0 & I\\R(1/\ell) & 0 &
0\endmatrix\right),$$  where $R(\theta)$ denote the standard
$2\times 2$ rotation matrix with rotation angle $2\pi \theta$, and
$I$ the $2\times 2$ identity matrix. The desired result follows.
\qed\enddemo

\vskip2mm

\demo{Proof of Proposition 1.4}

By Freedman [Fr], $S^5/T^1$ is homeomorphic to $\Bbb CP^2$.
Consider the induced $\Gamma$-action on $S^5/T^1$. Let
$\Gamma_0$ and $\Gamma^*_0$ be defined in the proof
of Lemma 1.7. Then $\Gamma _0\cong \Bbb Z_\ell$ is a subgroup
of $T^1$. By Theorem 1.5, the  $\Gamma^*$-action is conjugate
to a linear $\Gamma^*$-action on $\Bbb CP^2$ (and thus $
\Gamma^*$ is identified with a subgroup of $PSU(3)$, denoted by
$\Gamma^*_\ell$) by an equivariant homeomorphism $f:
(S^5/T^1,\Gamma^*)\to (\Bbb CP^2,\Gamma^*_\ell)$.

By Lemma 1.7, $\Gamma \cong \Gamma _\ell\subset SU(3)$, which acts
linearly on $S^5(1)$ that is the lifting of the $\Gamma^*_\ell$-action on
$\Bbb CP^2$ (but we should note that the $\Gamma _\ell$-action may
not be free a priorly). For the sake of convenience, let us
identify $\Gamma _\ell$ with $\Gamma$. It remains to prove that,
the free $\Gamma$-action on $S^5$ is conjugate to the linear
$\Gamma _\ell$-action on $S^5(1)$.

Consider the $(\Gamma, T^1)$- (resp. $(\Gamma _\ell, T^1)$)
principal bundle $S^5 \to S^5/T^1$.

By Theorem 4.5 (see Section 4) it suffices to prove that the induced principal
$T^1$-bundle $$T^1\to E\Gamma\times _\Gamma S^5\to E\Gamma \times
_\Gamma \Bbb CP^2 \tag 1.8$$ is equivalent to the corresponding
principal $T^1$-bundle of $(\Gamma _\ell, T^1)$-bundle on
$S^5(1)$. Clearly, $\pi _1(E\Gamma \times _\Gamma \Bbb
CP^2)=\Gamma$, and $\pi _2(E\Gamma \times _\Gamma \Bbb CP^2)\cong
\Bbb Z$. Therefore, $$H^2(E\Gamma\times _\Gamma\Bbb CP^2;\Bbb
Z)\cong \Bbb Z\oplus H_1(\Gamma)$$ where the free part may be
regarded as $\text{Hom}(H_2(E\Gamma\times _\Gamma\Bbb CP^2);\Bbb
Z)$ by the universal coefficients theorem.

Let $e_\Gamma$ denote the Euler class of the principal
$T^1$-bundle in (1.8). By the homotopy exact sequence one sees
that $e_\Gamma$ is a primitive element of $H^2(E\Gamma\times
_\Gamma\Bbb CP^2;\Bbb Z)$, i.e. modulo the torsion group
$H_1(\Gamma)$ it generates the group. Moreover, with the notions
in the proof of Lemma 1.7, $H_1(\Gamma)\cong \Bbb Z_{3\ell}$.

Let $\Bbb Z_{3\ell}:=\langle B\rangle \subset \Gamma$, which acts
on $\Bbb CP^2$ with three isolated fixed points. Let $[p]\in \Bbb
CP^2$ be such a fixed point. Consider the orbit $\Gamma [p]\subset
\Bbb CP^2$. The restriction of the fiber bundle (1.8) on $E\Gamma
\times _\Gamma \Gamma [p]=E\Bbb Z_{3\ell}\times _{\Bbb
Z_{3\ell}}[p]=B\Bbb Z_{3\ell}$ is equivalent to the principal
bundle
$$T^1\to E\Bbb Z_{3\ell}\times _{\Bbb Z_{3\ell}}T^1\to
E\Bbb Z_{3\ell}\times _{\Bbb Z_{3\ell}}[p] =B\Bbb Z_{3\ell}\tag
1.9$$ whose total space is homotopy equivalent to $T^1/{\Bbb
Z_{3\ell}} =T^1$. Therefore, the Euler class $e_\Gamma$, restricts
to the Euler class of (1.8), which is clearly a generator of
$H^2(E{\Bbb Z_{3\ell}}\times _{{\Bbb Z_{3\ell}}}[p])={\Bbb
Z_{3\ell}}$. On the other hand, the free part
$\text{Hom}(H_2(E\Gamma\times _\Gamma\Bbb CP^2);\Bbb Z)$ restricts
to zero in $H^2(E\Bbb Z_{3\ell}\times _{\Bbb Z_{3\ell}}[p])$.
Therefore, we may write $e_\Gamma $ as $$(1, \alpha )\in
\text{Hom}(H_2(E\Gamma\times _\Gamma\Bbb CP^2), \Bbb Z) \oplus
H_1(\Gamma) \cong \Bbb Z\oplus H_1(\Gamma)$$ where $\alpha \in
H_1(\Gamma)\cong \Bbb Z_{3\ell}$ is a generator.

Fix a generator $1\in H_1(\Gamma)\cong \Bbb Z_{3\ell}$. By [Wo]
Theorem 5.5.6 (cf. page 168, where $d=3$ for our case) there
always exists an automorphism $\psi=\psi _{1,t,u} :\Gamma\to
\Gamma$ (by sending $A$ to $A^t$, and $B$ to $B^u$) such that the
induced automorphism $[\psi]\in \text{Aut}(H_1(\Gamma))$ satisfies
that $[\psi ](\alpha )=\pm 1$ (depending $\mod (3)$ type of $t$).
Therefore, by composing the $\Gamma$-action with a suitable
automorphism $\psi$ of $\Gamma$, we may assume that
$e_\Gamma=(1,\pm 1)$.

The same goes through for the linear $\Gamma _\ell$-action on
$S^5(1)$. And we may assume its Euler class $e_{\Gamma _\ell}=(1,
\pm 1)\in H^2(E\Gamma _\ell \times \Bbb CP^2, \Bbb Z)$. It is easy
to see that, for the complex conjugated linear action of $\Gamma
_\ell$ on $S^5$, the Euler class is $(1, -1)$ (resp. $(1, 1)$), if
$e_{\Gamma _\ell}=(1, 1)$ (resp. resp. $(1, -1)$). Therefore the
$T^1$-principal bundle (1.8) is equivalent to the $T^1$-principal
bundle associated to some linear $\Gamma _\ell$-action on
$S^5(1)$. By Theorem 4.5 this implies that the $\Gamma$-action on
$S^5$ is heomorphically conjugate to a linear $\Gamma _\ell$-action
on $S^5(1)$. The desired result follows. \qed\enddemo

\vskip4mm

\subhead c. Proof of Theorem A by assuming Theorem E
\endsubhead

\vskip4mm

We need the following lemma to apply the criterion in Proposition
1.4.

\vskip2mm

\proclaim{Lemma 1.10}

Let $M$ be a closed $5$-manifold of positive curvature with
universal covering space $\tilde M\approx S^5$. Assume that $M$
admits a $\pi_1$-invariant isometric $T^1$-action. If $\pi_1(M)$
is not cyclic, then

\noindent (1.10.1) the $T^1$-action on $\tilde M$ is free;

\noindent (1.10.2) the $T^1$-action and the $\pi_1(M)$-action on $\tilde M$
commute.

\noindent (1.10.3) the induced $\pi_1(M)$-action on $\tilde M^*$ is
pseudo-free.
\endproclaim

\vskip2mm

\demo{Proof} If (1.10.1) does not hold, there is an isotropy group
$\Bbb Z_p\subset T^1$ whose fixed point set in $\tilde M$ is
either a circle or a homotopy $3$-sphere (see Theorem 4.3 in Section 4).
Clearly, $\pi _1(M)$ preserves the fixed point set $\tilde M
^{\Bbb Z_p}$. This trivially implies that $\pi _1(M)$ is cyclic if
$\text{dim}(M ^{\Bbb Z_p})=1$, a contradiction. If $\text{dim}(M
^{\Bbb Z_p})=3$ we get a totally geodesic $3$-manifold  $\tilde M
^{\Bbb Z_p}/\pi _1(M)$ in $M$. By Theorem 4.3 and [Ha], we see that
$\pi _1(M)$ is a spherical $3$-space group which
acts freely on $S^5$. We conclude once again that $\pi _1(M)$ is
cyclic, a contradiction.

If (1.10.2) is false, then there is an element $\gamma \in
\pi _1(M)$ such that the holonomy image $\rho (\gamma)=-1 \in
\text{Aut}(T^1)=\Bbb Z_2$. In other words, for any $t\in T^1$ and
$x\in \tilde M$ it holds that $t\gamma x=\gamma t^{-1}x$. Clearly
the induced action of $\gamma$ on $\tilde M^*=\Bbb CP^2$ has at
least a fixed point for a pure topological reason. Thus, $\gamma$
preserves a $T^1$-orbit, saying $T^1\cdot x$. Choose $t\in T^1$ so
that $t^2x=\gamma x$. Note that $\gamma tx=t^{-1}\gamma x=tx$.
This implies that $\gamma$ fixes the point $tx\in \tilde M$, a
contradiction to the fact that $\pi _1(M)$ acts freely on $\tilde
M$.

If (1.10.3) is false, then there is an element $\gamma \in \pi _1(M)$
whose action on $\tilde M^*=\Bbb CP^2$ has a $2$-dimensional
fixed point component $F_0$. Observe that, for any $x \in \tilde
M$ with image $x^*\in F_0$, there exists an element $t_{x^*}\in T^1$ (depending smoothly on $x^*\in F_0$) so that $t_{x^*}
\gamma (x)=x$. Clearly $t_{x^*}\in T^1$ has order $|\gamma |$
since $\Gamma$ acts freely on the orbit $T^1\cdot x$. Thus $t_{x^*}
$ is constant for all $x^*\in F_0$ since $F_0$ is connected.
This implies that $t_{x^*}\gamma$ contains the preimage of $F_0$
in its fixed point set, which has to be homeomorphic to $S^3$ by the
Smith theory and Hamilton's work once again. Therefore, $F_0$ must
be a $2$-sphere. The centralizer $C_{\langle \gamma \rangle}(\pi
_1(M)) $ of the subgroup $\langle \gamma \rangle \subset \pi
_1(M)$ acts freely on fixed point component ($\approx S^3$) of
$t_{x^*} \gamma$. For the same reasoning above we know that
$C_{\langle \gamma \rangle}(\pi _1(M)) $ is cyclic.

On the other hand, for any $\eta\in \pi_1(M)$, by Frankel's
theorem we then have that $\eta(F_0)\cap F_0\ne
\emptyset$, i.e. there is $x^*\in F_0$ such that $
\gamma(\eta(x^*))= \eta(x^*)$. This implies that $
\eta^{-1}\gamma \eta(x^*)=x^*$ and hence
$\eta^{-1}\gamma\eta$ and $\gamma$ both preserve the circle orbit
over $x^*$. Thus $\eta^{-1}\gamma\eta$ and $\gamma$  generate a
cyclic subgroup of $\pi _1(M)$. In particular,
 $\eta^{-1}\gamma\eta$ is in the centralizer $C_{\langle \gamma \rangle}(\pi _1(M)) $.
 Since $C_{\langle \gamma \rangle}(\pi _1(M)) $ is cyclic, clearly
$\eta^{-1}\gamma\eta$ (resp. $\gamma$) generates the unique cyclic
subgroup of order $|\gamma|$ in $C_{\langle \gamma \rangle}(\pi
_1(M)) $. Therefore, $\eta^{-1}\gamma\eta \in \langle\gamma
\rangle$. This proves that $\langle\gamma \rangle$ itself is a
normal subgroup in $\pi _1(M)$, and hence $\pi _1(M)$ acts freely
on the $S^3$ as above, which implies that $\pi _1(M)$ is cyclic. A
contradiction. \qed\enddemo

\vskip2mm

\demo{Proof of Theorem A by assuming Theorem E}

Consider $M$ as in Theorem A, whose Riemannian universal covering
space $\tilde M$ is homeomorphic to $S^5$ (the sphere theorem).
Let $S^5_\delta$ denote a sphere of constant curvature $\delta$.
By the volume comparison,
$$\text{vol}(M)=\frac{\text{vol}(\tilde M)}{|\pi_1(M)|}\le
\frac{\text{vol}(S^5_{1/4})}C<\epsilon$$ is small, and thus
without loss of generality we may assume $\tilde M$ admits a
$\pi_1(M)$-invariant isometric $T^k$-action (Theorems 1.1 and
1.2). By Theorem E, we may further assume that $k=1$. By Lemma
1.10, we can apply Proposition 1.4 to conclude the desired result.
\qed\enddemo

\vskip4mm

\head 2. Proof of Theorems B and C by Assuming Theorem E
\endhead

\vskip4mm

\subhead d. Proof of Theorem B by assuming Theorem E
\endsubhead

\vskip4mm

Consider $M$ in Theorem B. As in the proof of Theorem A, we may
assume that the volume of $M$ is small and thus $\tilde M$ admits a
$\pi_1(M)$-invariant isometric $T^k$-action. We will prove that
$k>1$ and then apply Theorem E.

Consider a $\pi_1$-invariant $T^1$-action on $\tilde M$. The
kernel of the holonomy representation, $\rho: \pi_1(M)\to
\text{Aut}(T^1)\cong \Bbb Z_2$, is a normal subgroup of index at
most two, and thus the $T^1$-action on $\tilde M$ descends to a
$T^1$-action on $\tilde M/\ker (\rho)$, which is either $M$ or a
double covering of $M$. In particular, a $T^1$-orbit on $M$ has a
tube on which there is induced $T^1$-action. We will call isotropy
groups of the local $T^1$-action isotropy groups of the
$\pi_1$-invariant $T^1$-action.

Let $\Cal Met$ denote the collection of isometric classes of
compact metric spaces. Equipped with the Gromov-Hausdorff distance
$d_{GH}$, $\Cal Met$ becomes a complete metric space. The Gromov
compactness ([Gr]) asserts that any sequence of closed $n$-manifolds
with Ricci curvature uniformly bounded from below and diameter
uniformly bounded above contains a convergent subsequence with
respect to $d_{GH}$ (note that the limit may not a be a Riemannian
manifold).

\vskip2mm

\proclaim{Lemma 2.1}

Assume a sequence of closed $n$-manifolds $M_i@>d_{GH}>>X$ such
that $|\text{sec}_{M_i}| \le 1$, where $X$ is a compact metric
space. If $\text{dim}(X)=(n-1)$, then there is a uniform upper
bound on the order of isotropy groups of the $\pi_1$-invariant
$T^1$-action on $M_i$.
\endproclaim

\vskip2mm

\demo{Proof} We argue by contradiction, assuming that $x_i\in M_i$
such that the isotropy group $T^1_{x_i}\cong \Bbb Z_{h_i}$ with
$h_i\to \infty$ (see (1.1.1)). Passing to a subsequence if
necessary, we may assume that $x_i\to x\in X$. Note that an open
neighborhood of $x$ is homeomorphic to a cone over the limit of
$S^\perp_{x_i}/T^1_{x_i}$, where $S^\perp_{x_i}$ is the unit
sphere in the normal space to $T^1(x_i)$, and $T^1_{x_i}$ acts on
$S^\perp_{x_i}$ via the isotropy representation. Because $h_i\to
\infty$, the limit of $S^\perp_{x_i}/T^1_{x_i}$ has dimension $\le
n-3$, and thus the cone has dimension $\le n-2$, a contradiction
to $\dim(X)=n-1$. \qed\enddemo

\vskip2mm

Consider an exceptional $T^1$-orbit $T^1(x)$ in $M$, with isotropy
group $\Bbb Z_h$. Then there is a lower bound on $h$ that is
related to the fundamental group. Let $\gamma$ denote the homotopy
class of $T^1(x)$ with order $r$, and let $\sigma$ be the homotopy
class of a principal $T^1$-orbit with order $s$. By analyzing the
covering map from a component in $\tilde M$ of the preimage of a
tube of $T^1(x)$ in $M$, one sees that $h\ge r/s$.

In the proof of Theorem B, we will use the following result on
fundamental groups of positively curved manifolds ([Ro3]). A
cyclic subgroup of $\pi_1(M)$ is called maximal, if it is not
properly contained in any cyclic subgroup of $\pi_1(M)$.

\vskip2mm

\proclaim{Theorem 2.2}

Let $M$ be a closed $n$-manifold of positive sectional curvature.
If $M$ admits a $\pi_1$-invariant isometric $T^k$-action, then any
maximal normal cyclic subgroup of $\pi_1(M)$ has index $\le w(n)$.
\endproclaim

\vskip2mm

\demo{Proof of Theorem B by assuming Theorem E}

By the volume comparison, $\text{vol}(M)=\text{vol}(\tilde
M)/|\pi_1(M)|\le \text{vol}(S^5_\delta)/w(\delta)$. We may assume
that $w(\delta)$ is large so that $\text{vol}(M)<\epsilon$. By
Theorems 1.1 and 1.2, without loss of generality we may assume
that $M$ admits a $\pi_1$-invariant isometric $T^k$-action. By
Theorem E, it suffices to show that $k>1$.

We argue by contradiction: assuming a sequence, $M_i$, satisfying
the above with $w_i(\delta)\to \infty$, and $k=1$. By the Gromov's
compactness, we may assume that $M_i@>d_{GH}>> X$. Because $k=1$,
it follows that $\dim(X)=4$, and thus any isotropy group of the
$\pi_1$-invariant $T^1$-action on $M_i$ has order $\le c$ (Lemma
2.1).

To get a contradiction, we will find an isotropy group of order
$>c$. Take any maximal normal cyclic subgroup
$H_i=\langle\gamma_i\rangle \subset \pi_1(M_i)$. By Theorem 2.2,
we obtain
$$\split w_i&\le [\pi_1(M_i):\text{cent}(\pi_1(M_i))]
\\&\le [\pi_1(M_i):H_i\cap \text{cent}
(\pi_1(M_i))]\\&=[\pi_1(M_i):H_i]\cdot [H_i:H_i\cap
\text{cent}(\pi_1(M_i))]\\&\le w(5)\cdot [H_i:H_i\cap
\text{cent}(\pi_1(M_i))].\endsplit$$ (note that the above implies
that $H_i$ is not trivial) Clearly, for $i$ large we may assume
that $[H_i: H_i\cap \text{cent}(\pi_1(M_i))]>c$. Let $\sigma_i$
denote the homotopy class of a principal $T^1$-orbit on $M_i$.
Then $\sigma_i$ is in the center of $\pi_1(M_i)$. Assume that
$\gamma_i$ preserves some $T^1$-orbit $T^1(\tilde x)$ (see Theorem 4.7 in
Section 4).
Because $\sigma$ preserves all $T^1$-orbits, $\gamma_i$ and
$\sigma_i$ generate a normal cyclic subgroup, and the maximality
of $H_i$ implies that $\sigma_i\in H_i$. Note that $\gamma_i$ is a
multiple of the homotopy class of the projection of $T^1(\tilde
x)$ in $M$. Then the isotropy group of $\pi(\tilde T^1(\tilde x)$
has order at least $[H_i:\langle \sigma_i\rangle ]\ge [H_i:H_i\cap
\text{cent}(\pi_1(M_i))]>c$, a contradiction. \qed\enddemo

\vskip4mm

\subhead e. Proof of Theorem C by assuming Theorem E
\endsubhead

\vskip4mm

Similar to the proof of Theorem B, Theorem C is a consequence of
the following proposition and Theorem E.

\vskip2mm

\proclaim{Proposition 2.3}

Let $M$ be a closed $n$-manifold of finite fundamental group satisfying
$$|\text{sec}_M|\le 1,\qquad \text{diam }(M)\le d,\qquad \frac{
\text{vol }(M)}{\max \{\text{injrad}(M,z)\}}<\epsilon_1(n,d).$$
Then $M$ admits a $\pi_1$-invariant isometric $T^k$-action with
$k>1$.
\endproclaim

\vskip2mm

For a motivation of Proposition 2.3, consider the metric product
of a unit sphere and a flat $\epsilon$-torus,
$M_\epsilon=S^n\times \epsilon^2T^k$. Then $
\frac{\text{vol
}(M_\epsilon)}{\max\text{injrad }(M_\epsilon)}\to 0$ (resp. is
proportional to $\text{vol }(S^n)$) if $k>1$ (resp. $k=1$), and
the $\pi_1$-invariant structure in Theorem 1.1 is the
multiplication on the $T^k$-factor.

\vskip2mm

\demo{Proof of Proposition 2.3}

We may assume that $\epsilon_1(n,d)$ is small so that $\text{vol
}(M)<\epsilon(n,d)$, and thus $M$ admits a $\pi_1$-invariant
almost isometric $T^k$-action (Theorem 1.1). Without the loss of
generality, we may assume that the metric is $T^k$-invariant
(Theorem 1.2).

We argue by contradiction; assuming a sequence, $M_i$, as in the
above such that $\text{vol}(M_i)/\max\{\text{injrad} (M_i,z)\}\to
0$ and $k=1$. Without loss of generality, we may assume that
$M_i@>d_{GH}>> X$. Let $x_i\in M_i$ such that $\text{injrad}
(M_i,x_i)=\max\{\text{injrad }(M_i,z)\}$. We claim that there is a
constant, $\eta>0$, such that, for all $i$, $T^1(x_i)$ is
contained in the $\eta$-tube $U_i$ of some $T^1$-orbit,
$T^1(y_i)$. Assuming the claim (whose proof is given at the end),
we will bound $\frac{\text{vol}(U_i)}{\text{max.}\text{injrad
}(M_i,x_i)}$ from below by a positive constant (depending on
$\eta$), a contradiction.

Because $T^1$ acts isometrically on $U_i$,
$$\text{vol}(U_i)=\text{length}(T^1(y_i))\cdot
\text{area}(D^\perp_i),\tag 2.3.1$$ where $D^\perp_i$ denotes a
normal slice of $U_i$. We shall bound $\text{area}(D^\perp_i)$
from below, and bound $\text{length}(T^1(y_i))$ in terms of
$\text{length}(T^1(x_i))$.

Let $T^1_{y_i}$ be the isotropy group of $T^1(y_i)$, and let $p_i:
\tilde U_i\to U_i$ denote the Riemannian $|T^1_{y_i}|$-covering
map, where $\tilde U_i=T^1\times D^\perp_i$, $U_i=T^1\times
_{T^1_{y_i}}D^\perp_i$ (the Slice lemma) and the lifting
$T^1$-action acts on $\tilde U_i$ by the rotation of the
$T^1$-factor. By the Gray-O"Neill Riemannian submersion formula,
the sectional curvature on $\tilde U_i/T^1$ is upper bounded by a
constant $c_1(n)$ (cf. (1.1.2)). By the volume comparison, we
conclude that $\text{vol} (D^\perp_i)=\text{area}(\tilde
U_i/T^1)\ge \text{vol} (B_\eta)$, where $B_\eta$ is a $\eta$-ball
in the $(n-1)$-space form of constant curvature $c_1(n)$. From
(2.3.1), we get
$$\text{vol}(U_i)\ge \text{length}(T^1(y_i))\cdot \text{vol}
(B_\eta).\tag 2.3.2$$

By (1.1.2) we may assume that the second fundamental group of all
$T^1$-orbits on $\tilde U_i$ are uniformly bounded by a constant
$c(n)\eta^{-1}$. Then we may assume a constant $c(n,r)$ such that
$\text{length}(T^1(\tilde x_i))\le c(n,\eta)\cdot
\text{length}(T^1(\tilde y_i))$, where $p_i(\tilde x_i)=x_i$ and
$p_i(\tilde y_i)=y_i$. Then
$$\text{length}(T^1(x_i))\le c(n,\eta)\cdot |T^1_{y_i}|\cdot
\text{length}(T^1(y_i)).\tag 2.3.3$$ Recall that the $T^1$-orbit
at any point represents all the collapsed directions of the metric
(cf. [CFG]). In particular, we may assume that
$$\text{injrad} (M_i,x_i)\le \frac 12\text{length} (T^1(x_i)).\tag
2.3.4$$

Using (2.3.2)-(2.3.4), we derive
$$\split \frac{\text{vol }(M_i)}{\text{injrad }
(M_i,x_i)}&\ge \frac{\text{vol }(U_i)}{\frac 12
\text{length}(T^1(x_i))}\\&\ge \frac{\text{length} (T^1(y_i))\cdot
\text{vol }(B_\eta)} {\frac 12 \text{length}
(T^1(x_i))}\\&\ge\frac{2\cdot \text{ vol }(B_\eta)}
{c(n,\eta)\cdot |T^1_{y_i}|}.\endsplit \tag 2.3.5$$ By Lemma 2.1,
we may assume that $|T^1_x|\le h(n,d)$ and thus see a
contradiction in (2.3.5) because the left hand side converges to
zero.

We now verify the claim. Recall that $M_i/T^1$ is homeomorphic and
$\epsilon_i$-isometric $X$, $\epsilon_i\to 0$, and the projection
of the singular set on $M_i$ into $X$ converging to the singular
set of $X$ (with respect to the Hausdorff distance). On the other
hand, $X$ is the metric quotient, $X=Y/O(n)$, where $Y$ is a
Riemannian manifold on which $O(n)$ acts isometrically (cf.
[CFG]). We can now pick up $\eta$ from the stratification
structure on $(Y,O(n)$ i.e. there are $O(n)$-invariant subsets,
$$Y=S_0\supset S_1\supset \cdots \supset S_r,\qquad \bar S_i=
\bigcup_{j\ge i} S_j$$ such that each component of $S_i$ has a
unique isotropy group, and $S_r$ is a closed totally geodesic
submanifold, where $\bar A$ denote the closure of a subset $A$.
We can choose a sequence of numbers,
$1>\eta_r>>\eta_{r-1}>> \cdots > \eta_1$ such that each $x\in
S_i=\bigcup_{j>i} T_{\eta_j}(S_j)$ satisfies that $O(n)(x)$ has a
$\eta_i$-tube. We then choose $\eta=\eta_1/2$. \qed\enddemo

\vskip2mm

\remark{Remark \rm 2.4}

Observe that the above proof goes through if one replaces the
assumption, ``$\frac{\text{vol}(M)}{\max \text{injrad}(M,z)}
<\epsilon_1(n,d)$'' by ``$\frac{\text{vol}(M)}{\text{injrad}(M)
}<\epsilon_1(n,d)$''. However, the later is not a valid assumption
since every spherical $5$-space form
satisfies $\frac{\text{vol}(S^5/\Gamma)}{\text{injrad}(S^5/\Gamma)}\ge
\frac{\text{vol}(S^5)}\pi$, which include the case `$k=1$'
(see the proof of Theorem C). On the other hand, given
$\epsilon>0$, there are
infinitely many spherical $5$-spaces satisfying $\frac{\text{vol}(S^5/\Gamma)}{\max
\text{injrad}(S^5/\Gamma,z)}<\epsilon$ (see Example 2.6).
\endremark

\vskip2mm

The following may be viewed as a converse to Proposition 2.3.

\vskip2mm

\proclaim{Lemma 2.5}

Let $M_i@>d_{GH}>>X$ such that $|\text{sec}_{M_i}|\le 1$ and
$\text{diam}(M_i)\le d$. If $\dim(X)\le n-2$, then
$\frac{\text{vol}(M_i)}{\max \{\text{injrad}(M_i,z)\}}\to 0$.
\endproclaim

\vskip2mm

\demo{Proof} We argue by contradiction; without loss of generality
we may assume that
$\frac{\text{vol}(M_i)}{\max\{\text{injrad}(M_i,z)\}}\ge c>0$ for all $i$.
By the Cheeger's lemma ([CE]), the ratio,
$$c\le \frac{\text{vol}(M_i)}{\text{injrad}(M)}\cdot
\frac{\text{injrad}(M_i)}{\max\{\text{injrad}(M_i,z)\}}=
\frac{\text{vol}(M_i)}{\max\{\text{injrad}(M,z)\}}\le c(n,d),$$
and thus $1\le
\frac{\max\{\text{injrad}(M_i,z)\}}{\text{injrad}(M_i)}\le
\frac{c(n,d)}c$.

Let $S_i$ denote the singular set of the $\pi_1(M_i)$-invariant
isometric $T^k$-action on $M_i$, and let $U_i$ denote the
$\epsilon$-tube of $S_i$. Then the orbit projection, $p_i:
M_i-U_i\to M_i^*-U_i^*$ ($x^*=p_i(x)$) is
a Riemannian submersion with fiber a flat manifold $F_i$. We may
choose $\epsilon$ small so that
$$\frac c2\le
\frac{\text{vol}(M_i-U_i)}{\max\{\text{injrad}(M_i,z)\}}
=\frac{\int_{M_i^*-
U_i^*} \text{vol}(p_i^{-1}(x^*))d\text{vol}^*}{\max\{\text{injrad}(M_i,z)\}}.$$ Because
$\text{diam}(p^{-1}_i(x^*))\to 0$ uniformly as $i\to \infty$,
again by the Cheeger's lemma we may assume that
$$\text{vol}(p_i^{-1}(x^*))\le \text{injrad}(p_i^{-1}(x^*))\epsilon_i,$$
where $\epsilon_i\to 0$ as $i\to \infty$.

Because the fiber $p^{-1}_i(x^*)$ points all collapsed
directions ([CFG]), we may assume that $\text{injrad}(M_i,x)\simeq
\text{injrad}(p_i^{-1}(x^*))$. Then
$$\split \frac c2&\le \frac{\text{vol}(M_i-U_i)}{\max\{\text{injrad}(M_i,z)\}}
\\& \le \int_{M_i^*-U_i^*}\frac{2\text{injrad}(M_i,x)}
{\max\{\text{injrad}(M_i,z)\}}\epsilon_id\text{vol}^*\\&\le
\frac{2c(n,d)}c\text{vol}(M_i^*-U_i^*)\epsilon_i\to
0,\endsplit$$ a contradiction. \qed\enddemo \vskip2mm

We now apply Lemma 2.5 to construct spherical $5$-space forms
satisfying Theorems B and C.

\vskip2mm

\example{Example 2.6}

We first construct lens spaces, $S^5/\Bbb Z_p$. Consider a
semi-free linear $T^2$-action on $S^5$. Let $T^1_i\subset T^2$
such that $\text{diam}(T^2/T_i^1)\le i^{-1}$ and $T^1_i$ has no
fixed point. We then choose a large prime $p_i$ so that $\Bbb
Z_{p_i} (\subset T^1_i)$ acts freely on $S^5$ and $\text{length}
(T^1/\Bbb Z_{p_i})\le i^{-1}$ (because the $T^1_i$-action has only
finitely many isotropy groups). Clearly, the Gromov-Hausdorff
distance, $d_{GH}(S^5/\Bbb Z_{p_i},S^5/T^2)\to 0$. By Lemma 2.5,
$S^5/\Bbb Z_{p_i}$ satisfies the conditions of Theorem C.

We now construct non-lens spaces. Let $A, B$ be the generators
of a spherical group as in Subsection b. Clearly, the center
of $\Gamma$ is generated by $\gamma^3_2$ and $[\Gamma:\langle \gamma^3_2\rangle ]\ge m$. Then $S^5/\Gamma$ satisfies conditions
of Theorems B and C when $m$ large.
\endexample

\vskip2mm

\example{Example 2.7}

We will construct examples showing that Theorems B and C are false
if relaxing ``$\delta>0$'' to ``$\delta\ge 0$''.

According to [Wo], there is a sequence of non-cyclic spherical
$3$-space groups, $\Gamma_i$, such that $S^3/\Gamma_i$ converges
to a closed interval $I$ and $[\Gamma_i:\text{cent}(\Gamma_i)]\to
\infty$. Then $M_i=S^2/\Gamma_i\times S^2$ converges to $I\times
S^2$ with $0\le \text{sec}_{M_i}\le 1$. However, $\Gamma_i$ cannot
act freely on $S^5$ (any finite group acting freely both on $S^3$
and $S^5$ must be cyclic, cf. [Bro]). Note that by Lemma 2.5,
$\frac{\text{vol}(M_i)}{\max\{\text{injrad}(M_i,z)\}}\to 0$.
\endexample

\vskip4mm

\head 3. Proof of Theorem D by Assuming Theorems E and F
\endhead

\vskip4mm

We first state a generalization of Theorem F. We call an
abelian subgroup $c_s$ of a finite group $\Gamma$ a
{\it semi-center}, if its centralizer has index at most two in
$\Gamma$ and if $|c_s|$ is `maximal' among all such abelian subgroups.
Obviously, a semi-center contains the center, and coincides with
the center when $|\Gamma|$ is odd. However, a semi-center may not
be unique when $|\Gamma|$ is even.

\vskip2mm

\proclaim{Theorem F'}

Let $M$ be a closed $5$-manifold of positive sectional curvature
which admits a $\pi_1$-invariant fixed point free isometric
$T^1$-action. Then the universal covering $\tilde M$ is
diffeomorphic to $S^5$, provided $\pi _1(M)$ satisfies
the following conditions:

\noindent (F1) $\pi_1(M)$ has a semi-center of index at
least $w>0$, a universal constant.

\noindent (F2) $\pi _1(M)$ does not contain any index $\le 2$
subgroup isomorphic to a spherical $3$-space group.
\endproclaim

\vskip2mm

We claim that Theorem F' implies Theorem F. Let $\pi_1(M)$
satisfies the assumptions of Theorem F. Then $\pi_1(M)$
is not cyclic. We shall check that $\pi_1(M)$ satisfies (F1) and (F2).
If $\pi_1(M)$ is isomorphic to a spherical
$5$-space group, then $\pi_1(M)$ contains an index $3$
normal cyclic subgroup and any semi-center
of $\pi_1(M)$ coincides with its center (see p.225, [Wo]).
Note that $\pi_1(M)$ contains no spherical $3$-space group of
index at most $2$; otherwise, the spherical $3$-space group
is cyclic because it also acts freely $S^5$ ([Bro]), and thus
$\pi_1(M)$ has a cyclic subgroup of index $2$, a contradiction.
If $\pi_1(M)$ has an odd order, then any semi-center has to
coincide with its center, and $\pi_1(M)$ has to be isomorphic
to a spherical $3$-space group which is cyclic because
its order is odd, a contradiction.

\vskip2mm

\demo{Proof of Theorem D by assuming Theorems E and F'}

Let $S^5_\delta$ denote a sphere of constant curvature $\delta$.
By the volume comparison, we may assume
$$\text{vol}(M)=\frac{\text{vol}(\tilde M)}{|\pi_1(M)|}\le
\frac{\text{vol}(S^5_\delta)}{w(\delta)}<\epsilon,$$
and thus we may assume, without loss of generality, that
$\tilde M$ admits a $\pi_1(M)$-invariant isometric
$T^k$-action (Theorems 1.1 and 1.2). If $k>1$, then by
Theorem E conclude the desired result. If $k=1$, then
by Theorem F' we see that the Riemannian universal covering of $M$
is diffeomorphic to a sphere, and thus by Proposition 1.4 and
Lemma 1.10 we
conclude the desired result.
\qed\enddemo

\vskip4mm

\head {\bf Part II. Proofs of Theorems E and F}
\endhead

\vskip4mm

\head 4. Preparations
\endhead

\vskip4mm

In this section, we supply materials that will be used in the
proofs of Theorems E and F' in the rest of this paper.

\vskip4mm

\subhead a. Fixed point set of abelian groups
\endsubhead

\vskip4mm

Consider a compact Lie group $G$ acting isometrically on a closed
manifold $M$. Let $M^G$ denote the set of $G$-fixed points.
Then each component of $M^G$ is a closed totally geodesic
submanifold. If $G=T^k$, then $F$ has even codimension. For a
generic compact Lie group, the topology of $M$ may not be well
related to the topology of $M^G$. However, the opposite
situation occurs when $G$ is abelian (cf. [Bre], p.163).

\vskip2mm

\proclaim{Theorem 4.1}

Let $M$ be a compact $\Bbb Z_h$-space. If $h=p$ is a prime, then
$$\chi(M)=\chi(M^{\Bbb Z_p}) \mod p.$$
If $\Bbb Z_p$ ($p$ is a prime) acts trivially on the
homology group $H_*(M;\Bbb Z)$, then
$$\chi (M)=\chi (M^{\Bbb Z_p}).$$
\endproclaim

\vskip2mm

The last assertion of the above lemma is from [Bre] III exercise
13 (p.169).

\vskip2mm

\proclaim{Theorem 4.2}

Let a compact abelian Lie group $G$ act effectively on a closed
manifold $M$, and let $N$ denote an invariant subset. Then
$$\operatorname{rank} (H^*(M^G,N^G;\ell))\le \operatorname
{rank}(H^*(M,N;\ell)),$$ where $G=T^k$ and $\ell=\Bbb Q$ or
$G=\Bbb Z^k_p$ and $\ell=\Bbb Z_p$.
\endproclaim

\vskip2mm

A consequence of Theorem 4.2 is

\vskip2mm

\proclaim{Theorem 4.3 (Smith)}

Let a torus $T^k$ act effectively on a closed manifold $M$. If
$M$ is a rational homology sphere, then $N^{T^k}$ is a rational
homology sphere.
\endproclaim

\vskip2mm

The $T^k$-action on a sphere without fixed points is well understood.

\vskip2mm

\proclaim{Theorem 4.4 ([Bre], P. 164)}

Let $M$ be a $n$-dimensional homology sphere and admit a $T^k$
action with no fixed point. If $H\subset T^k$ is a subtorus of
dimension $k-1$, let $r(H)$ denote that integer, for which
$M^H$ is a homology $r(H)$-sphere. Then with $H$ ranging over
all subtori of dimension $k-1$ and $r(H)\geqslant0$, we have
$$n+1=\sum_H(r(H)+1).$$
\endproclaim

\vskip2mm

\subhead b. A generalized Lashof-May-Segal theorem
\endsubhead

\vskip4mm

Let $G$ denote a compact Lie group ($G$ can be finite). For two
$G$-spaces, $Y$ and $Z$, a map,
$f: Y\to Z$, is called an {\it $G$-map} if $f(g\cdot y)=g\cdot f(y)$
for all $y\in Y$ and $g\in G$.

A {\it principal
$(G,T^k)$-bundle} is a principal $T^k$-bundle, $T^k\to E @>p>> Y$, such
that $E$ and $Y$ are $G$-spaces, $p$ is a $G$-map and the
$G$-action on $E$ preserves the structural group of the bundle.
Note that the $G$-action and $T^k$-action may not commute. Two
principal $(G, T^k)$-bundles are called {\it equivalent}, if there
is a $G$-equivariant bundle equivalent map.

Let $\Cal B(G,T^k)(B)$ be the set of equivalence classes of principal
$(G, T^k)$-bundles over $B$. When $G=\{1\}$, we will skip $G$ from the
notation.

Let $EG$ be the infinite join of $G$, a contractible free $G$-CW complex
(cf. [Hu]).  Put $B_G=EG\times _GB$ and $E_G=EG\times_GE$. There is a
natural transformation,
$$\Phi : \Cal B(G,T^k)(B)\to \Cal B(T^k)(B_G)$$
by sending a principal $(G,T^k)$-bundle $p: E\to B$ to the principal
$T^k$-bundle $p_G: E_G\to B_G$. The following theorem is a special
case of [FR2] Theorem 3.3 which generalizes the Lashof-May-Segal
theorem.

\vskip2mm

\proclaim{Theorem 4.5}

Let $B, G$ be as in the above. Then $\Phi: \Cal B(G,T^k)(B) \to \Cal
B(T^k)(B)$ is a bijection.
\endproclaim

\vskip2mm

Theorem 4.5 can be used in the following situation (see Section
7): Let $M$ be a closed manifold of finite fundamental group which
admits a pseudo-free $T^k$-action. Let $\tilde M_0=\tilde M-\Cal
S$, where $\Cal S$ is the union of singular orbits. Then
$\tilde M_0\to \tilde M_0^*$ is a
$(\pi_1(M),T^k)$-bundle.

\vskip4mm

\subhead c. Positive curvature and isometric torus actions
\endsubhead

\vskip4mm

In the rest of this section, we will consider an isometric
$T^k$-action on a closed manifold $M$ of positive sectional
curvature. As seen in Theorems 4.1--4.4, the topology of $M$ is
closely related to the singular structure and the orbit space of
the $T^k$-action. In the presence of a positive curvature, the
singular structure and the orbit space is very restricted.

A basic constraint on the singular structure is given by
following Berger's vanishing theorem ([Ro1], also [GS],
[Su]).

\vskip2mm

\proclaim{Theorem 4.6}

Let a torus $T^k$ act isometrically on a closed manifold $M$
of positive sectional curvature. Then there is a $T^k$-orbit
which is a circle. Moreover, the fixed point set is not
empty when $\dim(M)$ is even.
\endproclaim

\vskip2mm

Theorem 4.6 implies, via the isotropy representation at a circle
orbit, that large $k$ yields closed totally geodesic submanifolds
of small codimension.

In the study of the fundamental group of a positively curved
manifold on which $T^k$ acts isometrically, the following
result is a useful tool.

\vskip2mm

\proclaim{Theorem 4.7 ([Ro4])}

Let $M$ be a closed manifold of positive sectional curvature on
which $T^k$ acts isometrically, and let $\phi$ be an
isometry on $M$ which commutes with the $T^k$-action. Then $\phi$
preserves some $T^k$-orbit which is a circle.
\endproclaim

\vskip2mm

We now illustrate a situation where Theorem 4.7 may be applied.
Let $T^k$ act isometrically on a closed manifold $M$ of finite
fundamental group, and let $\pi: \tilde M\to M$ denote the
Riemannian universal covering. Let $p: \tilde M\to \tilde M^*$ be the orbit projection, where $\tilde T^k$ denotes the
covering torus of $T^k$ acting on $\tilde M$. For any $\gamma\ne
1\in \pi_1(M)$, because the $\gamma$-action commutes with the
$\tilde T^k$-action, $\gamma$ induces an isometry of $\tilde
M^*$, denoted by $\gamma^*$. Obviously,

\noindent (4.8.1) $\gamma^*$ is trivial if and only if $\gamma
\in H$, the subgroup generated by loops in a principal
$T^k$-orbit.

\noindent (4.8.2) $\gamma$ preserves an orbit, say $\tilde T^k
(\tilde x)$, if and only $\gamma^*$ fixes $x^*=p(\tilde x)$.

\noindent (4.8.3) If $k=1$ and $\gamma$ preserves $\tilde
T^1(\tilde x)$, then $T^1(\pi(\tilde x))$ is an exceptional orbit, whose
isotropy group contains a subgroup $\Bbb Z_h$ with $h$ the
order of $\gamma$.

\vskip2mm

\proclaim{Corollary 4.9}

Let $M$ be a closed manifold of positive sectional curvature. If
$M$ admits an isometric $T^k$-action, then the subgroup generated
by loops in a principal $T^k$-orbit is cyclic, say $\langle
\alpha\rangle$. Moreover, for any $\gamma\in \pi_1(M)$, $\alpha$
and $\gamma$ generate a cyclic subgroup.
\endproclaim

\vskip2mm

\demo{Proof} Let $H$ denote the subgroup generated by loops in a
principal $T^k$-orbit. Let $\tilde T^k$ denote the lifting
$T^k$ that acts on the Riemannian universal covering space $\tilde M$.
Then $H$ preserves all principal $\tilde T^k$-orbits and preserves
all $\tilde T^k$-orbit. Because there is a circle $\tilde T^k$-orbit,
$H=<\alpha>$ is cyclic. Because any $\gamma\in \pi_1(M)$ preserves some
circle orbits, $\alpha$ and $\gamma$ generates a cyclic group.
\qed\enddemo

\vskip2mm

The following connectedness theorem of Wilking provides a useful
tool to contract information on homotopy groups from the existence
of a closed totally geodesic submanifold of small codimension (cf. see [FMR]).

A map from $N$ to $M$ is called {\it $(i+1)$-connected}, if it induces
an isomorphism up to the $i$-th homotopy groups and a surjective
homomorphism on the $(i+1)$-th homotopy groups.

\vskip2mm

\proclaim{Theorem 4.10 ([Wi])}

Let $M$ be a closed $n$-manifold of positive sectional curvature,
and let $N$ be a closed totally geodesic submanifold of dimension
$m$. If there is a Lie group $G$ that acts isometrically on $M$
and fixes $N$ pointwisely, then the inclusion map is
$(2m-n+1+C(G))$-connected, where $C(G)$ is the dimension of a
principal orbit of $G$.
\endproclaim

\vskip4mm

\head 5. Proof of Theorem E for $k=3$
\endhead

\vskip4mm

Consider the case $k=3$ in Theorem E. By (1.3.1), we may assume
that the holonomy representation $\rho: \pi_1(M)\to\text{Aut}(T^3)
=\text{GL}(\Bbb Z,3)$ is not trivial and $\tilde M/\ker(\rho)$
is diffeomorphic to a lens space (in particular, $\ker (\rho)$
is cyclic). The goal of this section is to prove

\vskip2mm

\proclaim{Theorem 5.1}

Let $M$ be a closed $5$-manifold of positive sectional curvature.
If $M$ admits a $\pi_1$-invariant isometric $T^3$-action, then $M$
is homeomorphic to a spherical space form.
\endproclaim

\vskip2mm

By Proposition 1.4, the following two lemmas imply Theorem 5.1.

\vskip2mm

\proclaim{Lemma 5.2}

Let $M$ be a closed $5$-manifold of positive sectional curvature.
Suppose that $M$ admits a $\pi_1$-invariant isometric
$T^3$-action. Then $T^3$ has a circle subgroup $T^1$ which acts
freely on $\tilde M$ and which commutes with the
$\pi_1(M)$-action.
\endproclaim

\vskip2mm

\demo{Proof} From the above, $\tilde M\approx S^5$ and
$\text{ker}(\rho)$ is cyclic. We claim that $[\pi _1(M): \ker
(\rho)]=3$. It is easy to see that $\tilde M/T^3$ is homeomorphic
to a simplex $\triangle ^2$ as stratified set such that the three
vertices are the projection of isolated three circle orbits in
$\tilde M$ and the three edges are the projection of three
components of $T^2$-orbits (cf. [FR3]). If $\gamma\in \pi _1(M)$
acts trivially on the orbit space $\Delta ^2$, the holonomy $\rho
(\gamma)\in \text{Aut}(T^3)$ preserves all isotropy groups of the
$T^3$-action in $\tilde M$. In particular, $\rho (\gamma)$
preserves the three circle isotropy groups, $H_1\cap H_2, H_1\cap
H_3, H_2\cap H_3$, where $H_1, H_2, H_3$ are $2$-dimensional
isotropy groups. By the argument in the proof of (1.9.2) we know
that $\rho (\gamma)$ is the identity on all $H_1\cap H_2, H_1\cap
H_3, H_2\cap H_3$. Therefore, $\rho (\gamma)=I\in \text{GL}(\Bbb
Z, 3)$. On the other hand, an element $\gamma \in
\text{ker}(\rho)$ acts on $\Delta ^2$ preserving all vertices and
edges, and thus the identity in $\Delta ^2$. This clearly implies
that $\pi _1(M)/\text{ker}(\rho)\cong \Bbb Z_3$ acts effectively
on $\Delta ^2$ by rotating the three vertices.

Since $\rho (\pi _1(M))\subset \text{GL}(\Bbb Z, 3)$ is
non-trivial, by the above we know that $\rho (\pi _1(M))\cong \Bbb
Z_3$. Let $A$ denote a generator of $\rho (\pi _1(M))$. It is easy
to see that $A$ must have $1$ as an eigenvalue. Let $v$ be an
eigenvector of eigenvalue $1$, and let $T^s=\overline{\exp_etv}$
be the closed subgroup generated by the vector $v$. Then $T^s$
acts on $\tilde M$ commuting with the $\pi _1(M)$-action.

We claim that $T^s$ acts freely on $\tilde M$ and thus $s=1$. If
the claim is false, there is an isotropy group $H\ne 1\subset T^s$
with fixed point set $\tilde M ^H\ne\emptyset$. Note that $\tilde
M^H$ is either a circle or a totally geodesic three sphere
(Theorems 4.2 and 4.10). In either cases, $\tilde M^H$ is
invariant by the $T^3$-action and thus it contains one or two
circle orbits of the $T^3$-action. By the commutativity,
$\pi_1(M)$ also preserves $\tilde M^H$ and so $\pi _1(M)$
preserves one or two vertices of $\Delta ^2$. Therefore, $\pi
_1(M)$ has at least a fixed point among the three vertices of
$\Delta ^2$. A contradiction, since we have just shown $\pi_1(M)$
rotates the three vertices. \qed\enddemo

\vskip2mm

\proclaim{Lemma 5.3}

Let $T^1$ be as in Lemma 5.2. If $\rho : \pi _1(M)\to
\text{GL}(\Bbb Z,3)$ is non-trivial, then the induced
$\pi_1(M)$-action on $\tilde M^*\simeq \tilde M/T^1$ is pseudo-free.
\endproclaim

\vskip2mm

Note that Lemma 1.10 implies Lemma 5.3 if $\pi _1(M)$ is not
cyclic.

\demo{Proof of Lemma 5.3}

Because $\tilde M$ is diffeomorphic to $S^5$ on which $T^1$ acts freely,
$\tilde M^*\simeq \Bbb CP^2$ (cf. [Fr]). We argue by contradiction, assuming a $\gamma\in \pi_1(M)$ which has a $2$-dimensional fixed point set in $\tilde M^*$. Let $p: \tilde M\to \tilde
M/T^3 \approx \Delta ^2$ be the orbit projection. By the proof of
Lemma 5.2 we know that $\pi_1(M)/\text{ker}{ (\rho)}\cong \Bbb
Z_3$.

If $\rho (\gamma)$ is non-trivial, as we have seen in the proof of
Lemma 5.2, $\gamma$ rotates the vertices (and also edges) of
$\Delta ^2$ and hence it has only an isolated fixed point in the
interior of $\Delta ^2$. The preimage of this fixed point in
$\tilde M/T^1$ is a $2$-torus which must contain the fixed point
set of the $\gamma$-action on $\Bbb CP^2$. On the other hand, the
fixed point set of the $\gamma$-action on $\tilde M^*=\Bbb
CP^2$, which is a totally geodesic $2$-dimensional submanifold (a
sphere or $\Bbb RP^2$). A contradiction.

If $\rho (\gamma)$ is trivial, $\gamma$ acts trivially on $\Delta
^2$ by the proof of Lemma 5.2. For a $2$-dimensional connected
component $F_0$ in the fixed point set of $\gamma$-action in
$\tilde M^*$, there is an element $t_0\in T^1$ so that the fixed
point set of $t_0\gamma \in T^3\ltimes \pi _1(M)$ contains a
$3$-dimensional totally geodesic submanifold $F\subset \tilde M$
(cf. the proof of (1.10.3)). Since $\rho (\gamma )=I$, $F$ is
$T^3$-invariant with principal isotropy group a circle subgroup
$C\subset T^3$. Since $\text{ker} (\rho)\subset \pi _1(M)$ is
normal cyclic, $\langle \gamma \rangle $ is also a normal
subgroup. Thus $\pi _1(M)$ acts on the fixed point set $F_0$ and
so on $F$, i.e. $\alpha (F)=F$ for all $\alpha \in \pi _1(M)$. On
the other hand, for an element $\alpha \in \pi _1(M)$ with $\rho
(\alpha )$ non-trivial, the principal isotropy group of
$\alpha(F)$ is $\rho (\alpha) (C)$. Therefore $\rho (\alpha
)(C)=C$ and by the proof of Lemma 5.2, $C=T^1$. A contradiction, since $T^1$ acts freely on $\tilde M$. \qed\enddemo

\vskip4mm

\head 6. Proof of Theorem E at The Level of Fundamental Groups
\endhead

\vskip4mm

In this section, we will prove Theorem E at the level of
fundamental groups. The main result in this section is the
following:

\vskip2mm

\proclaim{Theorem 6.1}

Let $M$ be a closed $5$-manifold of positive sectional curvature.
If $M$ admits a $\pi_1$-invariant isometric $T^2$-action,
then the fundamental group of $M$ is isomorphic to that of a
spherical $5$-space form.
\endproclaim

\vskip2mm

By Theorem 5.1, it suffices to prove Theorem 6.1 for $k=2$.

\vskip2mm

\proclaim{Lemma 6.2}

A finite non-cyclic group $\Gamma$ is isomorphic to the
fundamental group of a spherical $5$-space form, if $\Gamma$
satisfies the following conditions:

\noindent (6.2.1) Every subgroup of order $3p$ is cyclic, for any
prime $p$.

\noindent (6.2.2) $\Gamma$ has a normal cyclic subgroup of index
$3$.
\endproclaim

\vskip2mm

\demo{Proof} By [Wo] Theorem 5.3.2 (in page 161), (6.2.1) and
(6.2.2) implies that every Sylow subgroup of $\Gamma$ is cyclic.
By [Wo] Theorem 5.4.1 (in page 163), $\Gamma$ is generated by
two elements $A$ and $B$ with relations
$$A^m=B^{n}=1, \qquad BAB^{-1}=A^r$$
where $((r-1)n, m)=1$ and $r^n\equiv 1(\text{mod }m)$. The order
$|\Gamma|=mn$. By (6.2.2), $n\equiv 0(\text{mod } 3)$, and $\{A,
B^3\}$ generates a normal cyclic subgroup of index $3$. Thus
$r^3\equiv 1(\text{mod }m)$. By [Wo] page 225, it only remains to
prove that $n\equiv 0(\text{mod }9)$. Suppose not, $(\frac n3,
3)=1$. For a prime factor $p|m$, $A^{\frac mp}$ is of order $p$.
By (6.2.1) the subgroup generated by $A^{\frac mp}$ and $B^{\frac
n3}$ is cyclic of order $3p$. On the other hand, this group is
a normal cyclic extension over a cyclic group satisfying the relation
$$B^{\frac n3}A^{\frac mp} B^{-\frac n3}=A^{{\frac mp}r^{\frac n3}}$$
This implies that $$r^{\frac n3} -1\equiv 0(\text{mod }p)\tag
6.2.3$$ Recall that $(r-1, p)=1$ and $r^3\equiv 1(\text{mod }p)$.
By considering $\frac n3(\text{mod }3)$ and (6.2.3) we get that
$r^2-1\equiv 0(\text{mod }p)$ and hence $r\equiv -1(\text{mod
}p)$. A contradiction, since $r^2+r+1\equiv 0(\text{mod }p)$.
 \qed\enddemo

\vskip2mm

In the proof of Theorem 6.1, we will establish (6.2.1) and
(6.2.2).

\vskip2mm

Recall that a $T^k$-action ($k>1$) is {\it pseudo-free}, if all
singular orbits are isolated.

\vskip2mm

\proclaim{Lemma 6.3}

Let $M$ be a closed $5$-manifold of positive sectional curvature
with a $\pi_1$-invariant isometric $T^2$-action. If the
$T^2$-action on $\tilde M$ is not pseudo-free, then $\pi_1(M)$ is
cyclic.
\endproclaim

\vskip2mm

\demo{Proof} By (1.3.2), $\tilde M$ is diffeomorphic to $S^5$.
Consider {\it all} circle subgroups of $T^2$ with nonempty fixed
point sets: By Theorem 4.4, there are three possibilities,
accordingly we divide the proof into three  cases:

Case (1). the fixed point set $\tilde M^{T^2}\ne \emptyset$;

Because $\pi _1(M)$ preserves the fixed point set, it suffices to
show that $\tilde M^{T^2}$ is a circle. By the Smith theory (cf.
Theorems 4.3) $\tilde M^{T^2}$ is a homology sphere. Since $T^2$
acts effectively on the normal space of $\tilde M^{T^2}$,
$\tilde M^{T^2}$ is a circle.

Case (2). There are two distinct circle subgroups: $T^1_1, T^1_2$
such that $\dim(\tilde M^{T^1_1})=3$ and $\tilde M^{T^1_2}$ is a
circle;

For any $\gamma \in \pi _1(M)$, observe that
$F(\rho(\gamma)(T^1_1), \tilde M)=\gamma(\tilde M^{T^1_1})$,
$\rho(\gamma)(T^1_1)=T^1_1$ where $\rho$ is the holonomy
representation (otherwise, there two distinct circle subgroups
with fixed point sets of dimension $3$). Consequently,
$\rho(\gamma)(T^1_2)=T^1_2$ and thus $\pi_1(M)$ preserves
$\tilde M^{T^1_2}$. Therefore, $\pi_1(M)$ is cyclic.

Case (3). There are three distinct circle subgroups with fixed
points set of dimension $1$, but there are non-trivial finite
isotropy groups.

Let $\Bbb Z_p\subset T^2$ ($p$ is a prime) be a finite isotropy
group. Note that  $\dim(\tilde M^{\Bbb Z_p})=3$. If
$\rho(\gamma)(\Bbb Z_p)=\Bbb Z_p$ for all $\gamma\in \pi_1(M)$,
then $\pi_1(M)$ preserves $\tilde M^{\Bbb Z_p}$. Because
$\tilde M^{\Bbb Z_p}$ contains exactly two of the three circle orbits in
Case (3), $\pi_1(M)$ must preserve the unique circle orbit outside
$\tilde M^{\Bbb Z_p}$ and thus $\pi_1(M)$ is cyclic.

If there is $\gamma\in \pi_1(M)$ such that $\rho(\gamma)(\Bbb
Z_p)\ne \Bbb Z_p$, then $\gamma(\tilde
M^{\Bbb Z_p})=\tilde M^{\rho(\gamma)(\Bbb Z_p)}$, and $F_0=\tilde M^{\Bbb Z_p}
\cap \tilde M^{\rho(\gamma)(\Bbb Z_p)}$ is a circle
which is the fixed point set of $\Bbb Z^2_p\subset T^2$ (Theorem
4.3). Because $\rho(\gamma)(\Bbb Z^2_p)=\Bbb Z^2_p$ for all
$\gamma\in \pi_1(M)$, $\pi_1(M)$ preserves $F_0$, and thus
$\pi_1(M)$ is cyclic. \qed\enddemo

\vskip2mm

\proclaim{Lemma 6.4}

Let $M$ be a closed $5$-manifold of positive sectional curvature.
Suppose that $M$ admits a $\pi_1$-invariant isometric $T^2$-action
with an empty fixed point set. If the $T^2$-action on $\tilde M$
is pseudo-free, then either $\pi_1(M)$ is cyclic or satisfies
(6.2.1) and (6.2.2).
\endproclaim

\vskip2mm

Combining Lemmas 6.2-6.4, we obtain the case of Theorem 6.1 for
$k=2$, and together with Theorem 5.1 we obtain Theorem 6.1.

For a pseudofree $T^2$-action on $\tilde M\approx S^5$, by Theorem
4.4 again there are exactly three isolated circle orbits.
Moreover, $\tilde M^*$ is a topological manifold and thus a
homotopy $3$-sphere because $\tilde M^*$ is simply connected.
Because $\gamma\in \pi_1(M)$ maps a singular orbit to
a singular orbit, we may view $\pi_1(M)$ acting on three points
(singular orbits) in $\tilde M^*$ by permutations. This defines
a homomorphism, $\phi: \pi_1(M)\to S_3$, the permutation group of
three letters. The kernel of $\phi$ is a normal subgroup acting
trivially on the three points (or equivalently, which preserving
every circle orbits). Thus $\text{ker}(\phi)$ is cyclic.

In the proof of Lemma 6.4 we need:

\vskip2mm

\proclaim{Lemma 6.5}

Let $M$ be as in Lemma 6.4. Then $\phi$ is trivial if and only if
the holonomy representation
 $\rho: \pi _1(M)\to \text{Aut}(T^2)$ is trivial.
\endproclaim

\demo{Proof} Let $H_i$, $i=1, 2, 3$, denote the three isotropy
groups of the isolated singular orbits of the pseudofree
$T^2$-action.  Note that, for any $\gamma \in \pi _1(M)$, and
$x\in \tilde M$ with isotropy group $I_x$, the isotropy group of
$\gamma (x)$, $I_{\gamma (x)}=\rho (\gamma ) (I_x)$. Therefore, if
$\rho$ is trivial, then $\pi _1(M)$ preserves the isotropy groups,
and so preserves every singular orbits, i.e., $\phi$ is trivial.

Conversely, if $\phi $ is trivial, $\pi _1(M)$ preserves the three
singular orbits. In particular, $\pi _1(M)$ is cyclic. Therefore,
$\rho (\gamma) (H_i)=H_i$ for any $\gamma \in \pi _1(M)$. It is
easy to see that $H_i$, $i=1, 2, 3$, generate $T^2$. Therefore, in
the Lie algebra of $T^2$, $\Bbb R^2$, the automorphism $\rho
(\gamma) \in \text{GL}(\Bbb Z, 2)$ has three different
eigenvectors whose eigenvalues are $1$ or $-1$. This implies that
$\rho (\gamma)=I$ or $-I$, where $I$ is the unit matrix. If $\rho
(\gamma)=-I$, i.e, for any $x\in \tilde M$, $\gamma
t^{-1}x=t\gamma x$. For a point $x$ in a singular orbit, it holds
$\gamma x=t^2x$ for some $t\in T^2$. Hence $\gamma t
x=t^{-1}\gamma x=tx$. This implies that $tx$ is a fixed point of
$\gamma$, a contradiction. \qed\enddemo

\vskip2mm

Our proof of Lemma 6.4 involves a homotopy invariant, `the first
$k$-invariant'. This invariant can be used to distinguish two
connected spaces whose first and second homotopy groups are the
same. Let's now briefly recall its definition ([Wh]).

Let $X$ be a connected space, and $K(\pi_i(X),\ell)$ denote the
Eilenberg-Maclane space. Corresponding to each map, $k_1: K
(\pi_1(X),1)\to K(\pi_2(X),3)$, there is a unique fibration,
$$\CD K(\pi_2(X),2)@ > >> E_{k_1}  \\
 @.  @ VVfV \\
 @.  K(\pi_1(X),1)@> k_1 >> K(\pi_2(X),3)
\endCD$$
with fiber $K(\pi_2(X),2)$. The total space $E_{k_1}$ is unique
and the classifying map $f$ has a lifting, $\tilde f: X\to
E_{k_1}$,
$$\CD @. E_{k_1} @.\\
 @.  @ VV V \\
 X @> f >>K(\pi_1(X),1)
\endCD$$
so that $\tilde f_*: \pi_i(X)\to \pi_i(E_{k_1})$ is an isomorphism
for $i= 1, 2$. The corresponding cohomology class $k_1\in
H^3(K(\pi _1(X),1); \pi _2(X))$ is called the first $k$-invariant
of $X$. Clearly, the first $k$-invariant is a homotopy invariant.

Let $L_\ell=S^3/\Bbb Z_\ell$ denote a lens space. It is well-known
that the punctured lens space has non-trivial first $k$-invariant,
i.e., for $p\in L_\ell$, $k_1(L_\ell-\{p\})\ne 0$ (cf. [EM]).

\demo{Proof of Lemma 6.4}

Let $\phi: \pi _1(M)\to S_3$ be the homomorphism as above. We
first claim that:

(6.4.1) the image of $\phi$ must be either $\{1\}$ or $\Bbb Z_3$.

Obviously, if $\text{Im}(\phi)\cong \{1\}$, then $\pi _1(M)$ fixes
all isolated circle orbits, thus $\pi _1(M)$ acts freely on every
circle orbit, and so $\pi _1(M)$ is cyclic. (6.2.1) clearly
follows by (6.4.1).

To prove (6.4.1) it suffices to rule out the case
$\text{Im}(\phi)\cong \Bbb Z_2$, since the case
$\text{Im}(\phi)=S_3$ may be ruled out by taking the pre-image
$\phi ^{-1}(\Bbb Z_2)$  of a order $2$ subgroup $\Bbb Z_2\subset
S_3$ instead of $\pi _1(M)$.

Assuming an element $\gamma\in \pi _1(M) $ so that $\phi (\gamma
)$ is of order $2$. By definition $\gamma$ preserves a unique
singular circle orbit, with isotropy group $H\cong S^1$, and
permutes the rest two singular circle orbits.
 Since $\tilde M^*$ is a homotopy $3$-sphere, the induced
action of $\gamma$ on
 $\tilde M^*$ has at least a fixed point which is not the
isolated singular points (the singular orbits). This implies that
$\gamma$ preserves a principal orbit $T^2\cdot x$ and acts freely
on. By Lemma 6.5 $\rho (\gamma)$ is also nonzero, of order $2$. It
is easy to show that the $\gamma$-action and the transitive
$T^2$-action on $T^2\cdot x$ do not commute. Therefore, the free
$\gamma$-action on $T^2\cdot x$ has a quotient space the Klein
bottle. This implies that, up to conjugation, the action of
$\gamma$ on the orbit is given by the composition of the
multiplication $\left [\matrix \alpha _1\\ \alpha _2
\endmatrix\right] \in T^2$ with the rotation $\left [\matrix 1 &
0\\ 0 & -1 \endmatrix\right]=\rho (\gamma)$ on $T^2$. Thus, for
$t_0=\left [\matrix \alpha ^{-1}_1\\ \alpha ^{-1}_2
\endmatrix\right]\in T^2$,
the fixed point set of $t_0\gamma$ on $T^2\cdot x$
contains two disjoint circles, $S^1\times \{\pm1\}$.
Therefore, the
fixed point set $F$ of $t_0\gamma$ on $\tilde M$ has
dimension $3$, a homology $3$-sphere, which
intersects with every principal orbit either empty, or
two disjoint circles.

Observe that $F$ projects to the fixed point set $\gamma$ on
$\tilde M^*$. Therefore, $\gamma$ acts on $\tilde M^*$ with
fixed point set a $2$-dimensional homology sphere. If $F$ does not
intersect with the circle orbit with isotropy group $H$ (preserved
by $\gamma$), the fixed point set of $\gamma$ on $\tilde M^*$ is
not connected, a contradiction. Since $H$ is preserved by $\rho
(\gamma)$, i.e. $\rho (\gamma )(H)=H$, $H$ is either $S^1\times
\{1\}$, or $\{1\}\times S^1$ in the standard coordinate for $T^2$.
Now the intersection of $F$ with the singular circle orbit
consists of either two points, or the whole singular orbit,
depending on whether the reduced automorphism $\rho (\gamma) \in
\text{Aut} (T^2/H)$ is trivial or not. In the former case, $H$
acts $F$ semifreely, with two isolated fixed points. A
contradiction, since $F$ is a homotopy $3$-sphere by Theorem 4.10,
because otherwise $H$ acts on two points punctured $3$-sphere
freely, absurd by Euler characteristic reasoning. For the latter
case, the quotient $F/H$ is a $2$-disk, with an action of $\{\pm
1\}\subset 1\times S^1$ freely in the interior of the disk. A
contradiction again by the Brower fixed point theorem.

It remains to prove (6.2.2). If $\pi_1(M)$ contains a non-cyclic
subgroup $G\subset \pi _1(M)$  of order $3q$ where $q$ is a prime,
by (6.2.1) $\phi (G)\cong \Bbb Z_3$ and thus $G$ contains a
subgroup $\Bbb Z_3$ acting pseudo-freely on $\tilde M^*$
(denoted by $\Sigma$), a homotopy $3$-sphere. Let $\tilde M_0$
denote the complement of the three isolated circle orbits on
$\tilde M$. Consider the
$T^2$-bundle,
$$T^2\to \tilde M_0/\Bbb Z_3\to \tilde M^*_0/\Bbb Z_3,$$
and the classifying map $f: \tilde M^*_0/\Bbb Z_3\to B(T^2\ltimes \Bbb
Z_3)$ of the principal $T^2\ltimes \Bbb Z_3$ bundle $\tilde M_0\to
\tilde M^*_0/\Bbb Z_3$. Because $\tilde M_0$ is $2$-connected, $f$ is
a $3$-equivalence. This implies that the first $k$-invariant of
$\tilde M^*_0/\Bbb Z_3$ is zero, a contradiction, because
$\tilde M^*_0/\Bbb Z_3$ is homeomorphic to the punctured lens space
$L_3-\{p\}$ which has a non-zero first $k$-invariant.
\qed\enddemo

\vskip4mm

\head 7. Proof of Theorem E for Pseudo-free $T^2$-actions
\endhead

\vskip4mm

By Theorem 5.1, the proof of Theorem E reduces to the
case of a $\pi_1$-invariant isometric $T^2$-action. The
goal of this section is to prove the case of a pseudo-free
$\pi_1$-invariant $T^2$-action (Theorem 7.1). In a non
pseudo-free situation, a totally geodesic submanifold of
codimension $2$ is present that requires a different
argument (see Section 8).

A crucial information in the proof of pseudo-free case is
that the orbit space $\tilde M^*$ is a homotopy $3$-sphere.
A serious problem is that the Poincar\'e conjecture is open
and we can not conclude that it is homeomorphic to $S^3$.
This will be solved because everything can go through by employing
$s$-cobordism theory in differential topology, and using the
well-known $s$-cobordism theorem (due to Smale) in dimension $5$.

In this section, for the sake of a simple exposition, we will
present a proof by assuming the Poincar\'e conjecture and leave
the proof in the general case to the Appendix.

\vskip2mm

\proclaim{Theorem 7.1}

Let the assumptions be as in Theorem E. If $k=2$ and the
$T^2$-action is pseudo-free, then $M$ is homeomorphic to a
spherical space form.
\endproclaim

\vskip2mm

Before starting a proof, it may be helpful to look at a
linear model in Theorem 7.1. Consider a linear pseudo-free
$T^2$-action as in {\bf b.} of Section 1, and
a commuting linear $\Bbb Z_{k\ell}$-action on
$S^5(1)$. Then $\Bbb Z_{k\ell}$ also acts on the orbit space $S^5(1)/T^2=S^3$. Let $\Bbb Z_k$
denote the principal isotropy group of the reduced linear action
on $S^3$. We will call the $\Bbb Z_k$-action {\it along the
$T^2$-orbits}. The reduced $\Bbb Z_\ell=\Bbb Z_{k\ell}/\Bbb Z_k$
action on $S^3$ has a fixed point set $S^1$, which spans a plane
of $\Bbb R^4$. The condition of the linear $\Bbb Z_k$ action
(defined by multiplying $(e^{\frac ak\pi i},e^{\frac bk\pi i},
e^{\frac ck\pi i}$)) along the $T^2$-orbits can be written as $
{a+b+c}=0(\text{mod } 2k)$. Note that in the above model, the action
by $A$ is always along the $T^2$-orbits.

In the following context we will continue to use $\Gamma_\ell$ (resp.
$T^2_\ell$) to mean a linear $\Gamma$-action ($T^2$-action) on $S^5(1)$
so that it extends to a linear action of $T_\ell\ltimes \Gamma$.

Consider a pseudofree linear $T^2$-action on $S^5$ defined in the
above. Let $S^5_0$ denote the complement of small open tubes
($\simeq D^4\times T^1$) around the three isolated circle orbits
on $S^5$. Let $M$ be as in Theorem 7.1. By Theorem 4.4 the
pseudofree $T^2$-action on $\tilde M\approx S^5$ has exactly three
isolated circle orbits. Let $\tilde M_0$ denote the complement of
small open tubes of the three isolated circle orbits. By Theorem
6.1 $\Gamma=\pi_1(M)$ is a spherical $5$-space group. Our main
effort is to show that $M_0=\tilde M/\Gamma$ is homeomorphic to
$S^5_0/\Gamma_0$. By [SW] the gluing of a handle $D^4\times T^1$
is unique up to homeomorphism, and therefore $M$ is homeomorphic
to $S^5/\Gamma _0$.

Let  $\rho :\Gamma\to \text{Aut} (T^2)$ denote the holonomy
representation of the $\pi_1$-invariant action. By Lemma 6.6 we
know that $\text{ker}(\rho)$ is cyclic, and the image
$\rho(\Gamma)$ is either trivial or isomorphic to $\Bbb Z_3$. Let
$\Bbb Z_k$ denote the principal isotropy group of the reduced
$\Gamma$-action on $\tilde M^*$. By definition one sees
that $\Bbb Z_k$ acts on $\tilde M$ through the $T^2$-orbits.

Now let us consider the principal $T^2$-bundle $T^2\to \tilde
M_0\to \tilde M^*_0$. Assuming the Poincare
conjecture, $\tilde M^*_0=S^3_0$ is the complement of $S^3$ by
removing three small $3$-disks. Note that $\Gamma$ acts on this
principal bundle, and moreover, the sub-action of
$\text{ker}(\rho)$ commutes with the $T^2$-action. The following
lemma is immediate:

\vskip 2mm

\proclaim{Lemma 7.2 }

The principal $T^2$-bundle is unique up to
weak equivalence. Therefore, every pseudofree $T^2$-action on
$S^5$ is conjugate to a linear $T^2$-action.
\endproclaim

\vskip2mm

\demo{Proof} Note that the bundle is uniquely determined by its
Euler class, an element in $H^2(\tilde M^*_0; \Bbb Z^2)\cong
\text{Hom}(\Bbb Z^2,\Bbb Z^2)$. Considering the Euler class as a
$2\times 2$ matrix (given by its classifying map to $BT^2$), its
determinant is $\pm 1$ since the total space $\tilde M_0$ is
$2$-connected by the transversality theorem. Therefore, up to the
left action by $GL(\Bbb Z, 2)$, i.e. up to an automorphism of
$T^2$, the bundle is unique. The desired result follows.
\qed\enddemo

\vskip2mm

By Lemma 7.2 the sub-action by $\Bbb Z_k$ on $\tilde M\approx S^5$
is conjugately linear. Therefore $\tilde M/\Bbb Z_k$ is
diffeomorphic to a lens space $S^5/\Bbb Z_k$. Of course one should
note that there are possibly many different ways to embed $\Bbb
Z_k$ in $T^2$ which acts freely on $S^5$, and consequently the
lens space may not be unique.

Let us consider the reduced principal $T^2$-bundle $T^2/\Bbb
Z_k\to \tilde M_0/\Bbb Z_k \to \tilde M^*_0$, regarded as a
$\Gamma/\Bbb Z_k$-equivariant bundle.

\vskip 2mm

\proclaim{Lemma 7.3}

If the $\Gamma$-action and $T^2$-action on
$\tilde M$ commute, then the above $\Gamma$-equivariant principal
$T^2$-bundle is $\Gamma/\Bbb Z_k$-equivariantly equivalent to a linear
$T^2$-bundle $T^2/\Bbb Z_k\to S^5_0/\Bbb Z_k\to S^3_0$.
\endproclaim

\vskip2mm

\demo{Proof} Note that $\Gamma$ is cyclic. Let us write $\Gamma
=\Bbb Z_{k\ell}$ where $\Bbb Z_k$ is as above. The effective
action on $\Bbb Z_\ell$ on $\Sigma \approx S^3$ has fixed point
e.g., the three isolated marked points representing the three
singular orbits. By a deep theorem of [BLP] this action of $\Bbb
Z_\ell$ on $\Sigma$ is conjugate to a linear action on $S^3$.
Therefore, by Theorem 4.5 it suffices to prove that the associated
principal $T^2$-bundle
$$
T^2/\Bbb Z_k\to E\Bbb Z_\ell\times _{\Bbb Z_\ell}\tilde M_0/\Bbb
Z_k\to E\Bbb Z_\ell\times _{\Bbb Z_\ell} \tilde M^*_0
$$
is unique up to weak equivalence. We need only to show its Euler
class $e(\gamma) \in H^2(E\Bbb Z_\ell\times _{\Bbb Z_\ell}\Sigma
_0;\Bbb Z^2)$ can be realized by the Euler class of a linear
$\Gamma _0$-equivariant principal $T^2$-bundle over $S^3_0$ with
total space $S^5_0/\Bbb Z_k$, where $\Bbb Z_k$ acts linearly on
$S^5$ along the $T^2$-orbits. By some standard calculation we get
that $H^2(E\Bbb Z_\ell\times _{\Bbb Z_\ell}\tilde M^*_0;\Bbb
Z^2)\cong \text{Hom}(\Bbb Z^2, \Bbb Z^2)\oplus \Bbb Z_\ell ^2$. By
restricting the bundle to $E\Bbb Z_\ell \times _{\Bbb Z_\ell} [p]$
where $[p]\in \tilde M^*_0$ is a fixed point of the $\Bbb
Z_\ell$-action, one sees that $e(\gamma)$ restricting to
$H^2(B\Bbb Z_\ell;\Bbb Z^2)\cong \Bbb Z_\ell^2$ is an
element of order $\ell$.

By comparing with the linear model discussed at the beginning of
this section, it is straightforward to check that every pair $(a,
b)\in \Bbb Z_\ell ^2$ generating an order $\ell$ element can be
realized as the torsion component of the Euler class of a linear
$\Bbb Z_\ell$-equivariant principal $T^2$-bundle on $S^3_0$ with
total space $\tilde M_0/\Bbb Z_{k\ell}$. The torsion free part of
$e(\gamma)$ is uniquely determined by its lifting to $H^2(E\Bbb
Z_\ell\times \tilde M^*_0;\Bbb Z^2)\cong \text{Hom}(\Bbb Z^2, \Bbb
Z^2)$, which is the Euler class of the forgetful principal
$T^2$-bundle on $\tilde M^*_0$, regarded as a non-equivariant bundle.
By Lemma 7.2 this Euler class is uniquely determined by the total
space $\tilde M_0/\Bbb Z_k$, or equivalently, by the conjugacy
class of the embedding of $\Bbb Z_k$ in $T^2$. This proves the
desired result.
\qed\enddemo

\vskip2mm

Next let us consider the case where the holonomy $\rho: \Gamma \to
\text{Aut}(T^2)$ is non-trivial.

\vskip2mm

\proclaim{Lemma 7.4}

Let $M$ be as in Theorem 7.1. If the holonomy
$\rho:\Gamma \to \text{Aut}(T^2)$ is non-trivial, then the
$\Gamma$-equivariant principal $T^2$-bundle $ \tilde M_0/\Bbb
Z_k\to \tilde M^*_0$ is $\Gamma/\Bbb Z_k$-equivariantly equivalent to a linear
$T^2$-bundle $T^2/\Bbb Z_k\to S^5_0/\Bbb Z_k\to S^3_0$.
\endproclaim

\vskip2mm

We first need some preparation. By Lemma 6.5 the image
$\rho(\Gamma)\cong \Bbb
Z_3\subset \text{GL}(2,\Bbb Z)=\text{Aut}(T^2)$. Recall that
$\text{SL}(2,\Bbb Z)\cong \Bbb Z_4*_{\Bbb Z_2}\Bbb Z_6$ has a
subgroup of order $3$, unique up to conjugation, which is
generated by
$$\left [\matrix 0 & -1 \\ 1 & -1 \endmatrix\right].$$
By Lemma 6.5 there exists an element $\gamma \in \Gamma$ such that
$\phi (\gamma)$ has order $3$.

Let $\Bbb Z_k\subset \Gamma$ denote
the principal isotropy group of the induced action on $\Sigma$.
The quotient group $\Gamma/\Bbb Z_k$ acts effectively on $\Sigma
\approx S^3$, and which is not free, unless $\Gamma /\Bbb Z_k\cong
\Bbb Z_3$. By Theorem 6.1 it is easy to see that $\Gamma/\Bbb
Z_k\cong \Bbb Z_3$ only if $\Gamma$ is cyclic. Therefore, by [BLP]
once again the reduced action on $\Gamma/\Bbb Z_k$ on $\Sigma
\approx S^3$ is conjugate to a linear action on $S^3$, unless
$\Gamma/\Bbb Z_k=\Bbb Z_3$. In the latter case (where $\Gamma/\Bbb
Z_k=\Bbb Z_3$) we may replace the "conjugation" by "s-cobordism"
(cf. Appendix), and everything goes through. For the sake of
simplicity we now assume that $\Gamma/\Bbb Z_k$ acts linearly on
$\Sigma\approx S^3$. Therefore, $\Gamma/\Bbb Z_k$ is a subgroup of
$SO(4)$.

\vskip2mm

\proclaim{Sublemma 7.5}

$\Gamma/\Bbb Z_k$ is cyclic.
\endproclaim
\demo{Proof} We need only to consider the case where $\Gamma$ is
not cyclic. By Theorem 6.1 $\Gamma=\{A, B: A^m=B^n=1,
BAB^{-1}=A^r\}$ is a spherical $5$-space group where $n\equiv
0(\text{mod }9)$, $(n(r-1), m)=1$, and $r^3\equiv 1(\text{mod
}m)$. Observe that $m$ must be odd. Hence, any possible
$2$-subgroup of $\Gamma$ is cyclic, saying $\Gamma _2$, such that
$\Gamma =\Gamma _2\times \Gamma /\Gamma _2$, where $\Gamma /\Gamma
_2$ is again a spherical $5$-space group of odd order. Since every
odd order subgroup of $SO(4)$ is also a subgroup of
$\text{Spin}(4)=S^3\times S^3$, by the classification of finite
subgroups of $S^3$ one concludes that a odd order subgroup of
$S^3\times S^3$ is abelian of rank at most $2$. Applying this we
know that $\Gamma /\Bbb Z_k$ is abelian. By the presentation of
$\Gamma$ above this implies that $\Gamma/\Bbb Z_k$ is
cyclic.
\qed\enddemo

\vskip2mm

Let us write $\Gamma /\Bbb Z_k=\Bbb Z_{3\ell}$.

\vskip2mm

\proclaim{Sublemma 7.6}

If $k\ne 1$, then $\Gamma$ is not cyclic and has a normal cyclic
subgroup whose quotient group is cyclic, i.e.,  $1\to \Bbb Z_m\to
\Gamma\to \Bbb Z_n\to 1$.
\endproclaim

\vskip2mm

\demo{Proof} Recall that $\Bbb Z_k$-acts in $\tilde M$ as a
sub-action of the pseudo-free $T^2$-action. Since the $T^2$-action
is $\Gamma$-invariant, for any $\alpha \in \Gamma$, it holds that
$\alpha \beta \alpha ^{-1}=\rho (\alpha) (\beta)$, where $\beta
\in \Bbb Z_k\subset T^2$ (Notice $\Bbb Z_k$ is invariant by $\rho
(\alpha)$). If $\Gamma$ is abelian, then $\rho (\alpha )(\beta
)=\beta$ even if $\rho (\alpha )=\left [\matrix 0 & -1 \\ 1 & -1
\endmatrix\right]\in \text{Aut}(T^2)$. A contradiction, since $1$ is not an eigenvalue of
the matrix. It is clear that $\Gamma$ satisfies the desired property.\qed\enddemo

\vskip2mm

\demo{Proof of Lemma 7.4}

By the discussion at the beginning of this section, there is a
linear (and free) $\Gamma$-action on $S^5(1)$ so that the matrix
$A$ has order $k$ which acts along the $T^2$-orbits (with the
pseudofree $T^2$-action as in the beginning), and $B$ has order
$3\ell$, where $\ell$ is divisible by $3$. It is easy to see that
the linear $T^2$-action is $\Gamma$-invariant. Since the linear
$\Bbb Z_{3\ell}$-action on $S^3$ has no fixed point but the
subgroup $\Bbb Z_\ell$ has fixed point set a trivial knot in
$S^3$, it is easy to see that such a linear action is unique up to
conjugation.

By Theorem 4.5 the affine $(\Bbb Z_{3\ell },T^2)$-bundle $\tilde
M_0/\Bbb Z_k\to \tilde M^*_0$, is uniquely determined by the
associated affine $T^2$-bundle with the above holonomy $\rho :\Bbb
Z_{3\ell}\to \text{GL}(\Bbb Z, 2)$, the latter is characterized by
its Euler class, an element in the local cohomology group
$H^2(E\Bbb Z_{3\ell}\times _{\Bbb Z_{3\ell}}\tilde M^*_0;\Bbb
Z^2_\rho)$. To prove lemma 7.4 we only need to verify that the
Euler class of the affine bundle is uniquely determined by Euler
class of the principal $T^2$-bundle $\tilde M_0/\Bbb Z_k\to \Sigma
_0$, that implies $ \tilde M_0/\Bbb Z_k\to \tilde M^*_0$ is
$\Gamma$-equivariantly equivalent to the data in the linear model
case.

By the short exact sequence $1\to \Bbb Z^2_\rho \to \Bbb Z[\Bbb
Z_3] \to \Bbb Z\to 1$ we can calculate the local cohomology group
$$1\to H^2(E\Bbb Z_{3\ell} \times_{\Bbb
Z_{3\ell}} \tilde M^*_0; \Bbb Z^2_\rho)\to H^2(E\Bbb Z_\ell
\times_{\Bbb Z_\ell} \tilde M^*_0; \Bbb Z)\to  H^2(E\Bbb Z_{3\ell}
\times_{\Bbb Z_{3\ell}} \tilde M^*_0; \Bbb Z)$$ where the middle
space $E\Bbb Z_\ell \times_{\Bbb Z_\ell} \tilde M^*_0$ is the three
fold covering of $E\Bbb Z_{3\ell} \times_{\Bbb Z_{3\ell}} \Sigma
_0$. In the above exact sequence, the middle term is isomorphic to
$\Bbb Z^2 \oplus \Bbb Z_\ell$, and the last term is isomorphic to
$\Bbb Z_{3\ell}$. By the universal coefficients theorem it is
readily to see that the torsion part of the middle term goes
injectively into the last term. Therefore, the local cohomology
group $H^2(E\Bbb Z_{3\ell}\times _{\Bbb Z_{3\ell}}\tilde M^*_0;\Bbb
Z^2_\rho)$ is torsion free of rank $2$. Hence, by the universal
coefficients theorem again, $H^2(E\Bbb Z_{3\ell}\times _{\Bbb
Z_{3\ell}}\tilde M^*_0;\Bbb Z^2_\rho)$ is given by $\text{Hom}^\rho
(H_2(\tilde M^*_0),\Bbb Z^2)\cong\text{Hom}^\rho (\Bbb Z^2,\Bbb Z^2)$,
where $\text{Hom}^\rho$ denotes the $\rho$-invariant
homomorphisms. Therefore, the forgetful homomorphism
$\text{Hom}^\rho (\Bbb Z^2,\Bbb Z^2)\to \text{Hom}(\Bbb Z^2,\Bbb
Z^2)$ sends the Euler class of the $\Bbb Z_{3\ell}$-equivariant
$T^2$-bundle on $\tilde M^*_0$ to the Euler class of the principal
$T^2$-bundle (forgetting the $\Bbb Z_{3\ell}$-action), which is
clearly injective. The desired result follows.
\qed\enddemo

\vskip2mm

\remark{Remark \rm 7.7}

By the proof of Lemma 7.4 and the
comments after Lemma 7.2, we know that the
$\Gamma$-equivariant principal $T^2$-bundle in
 Lemma 7.4 is uniquely determined by the total space $\tilde
 M_0/\Bbb Z_k$ and so by the lens space $\tilde M/\Bbb Z_k$,
 which depends only on the embedding of $\Bbb Z_k$ in $T^2$.
\endremark

\vskip2mm

\demo{Proof of Theorem 7.1}

Let $M$ be as in Theorem 7.1. By Theorem 6.1 $\pi _1(M)=\Gamma $
is a spherical $5$-space group. By Lemmas 7.3 and 7.4 we know that
$\tilde M_0/\Gamma :=M_0$ is diffeomorphic to $S^5_0/\Gamma $.
Since $M$ is obtained by gluing three handles $S^1\times D^4$
along the boundary components $S^1\times S^3$. Because every self
diffeomorphism of $S^1\times S^3$ extends to a self homeomorphism
of $S^1\times D^4$ (cf. [SW]), the homeomorphism type does not
depend on the gluing. Therefore, $M$ is homeomorphic to
$S^5(1)/\Gamma $. The desired result follows. \qed\enddemo

\vskip4mm

\head 8. Completion of The proof of Theorem E
\endhead

\vskip4mm

After the works in Sections 6-7, we are ready to finish the
remaining case in the proof of Theorem E, i.e., the case that a
$\pi _1$-invariant $T^2$ action on $M$ is not pseudofree.

\vskip2mm

\demo{Proof of Theorem E}

First, by Theorems 7.1 we only need to consider a non-pseudo-free
$T^2$-action, i.e. the $T^2$-action has a non-empty fixed point, a
$3$-dimensional stratum with circle isotropy group, or a
nontrivial finite isotropy group but without fixed point. In all
cases, $\pi_1(M):=\Gamma$ is cyclic (Lemmas 6.3 and 6.4). Recall
that $\tilde M=S^5$ with an action of $T^2\ltimes_\rho \Gamma$.

Case 1. The $T^2$-action has a non-empty fixed point set.

Note that the fixed point set must be a circle (Theorem 4.2). By
local isotropy representation of $T^2$ at the fixed point set,
there are two circle isotropy groups with three dimensional fixed
point sets, two totally geodesic $S^3$. Observe that the
$T^2$-action on $\tilde M$ is free outside the union of the two
$3$-dimensional strata, and the quotient space $\tilde M^*$ is
homeomorphic to the $3$-ball $D^3$, whose boundary $S^2=D^2_+\cup
D^2_{-}$, where $D^2_{\pm}$ is the image of the two
$3$-dimensional strata and $D^2_+\cap D^2_{-}$ is the image of the
fixed point set of $T^2$. We claim that {\it the $T^2$-action and
$\Gamma$-action commute} (equivalently, $\rho $ is trivial). By
now it is easy to see that $M$ is diffeomorphic to a lens space
(cf. [GS]).

Identify $\tilde M^*$ with $D^3$. Observe that $\Gamma$ acts
isometrically on $\tilde M^*$ and preserves the boundary
$\partial (\tilde M^*)=\partial D^3$. If the commutativity
fails, $\Gamma$ acts non-trivially on $\tilde M^*$. By the
well-known Brouwer fixed point theorem, $\Gamma$ has at least a
fixed point in the interior of $D^3$, which represents a principal
orbit, say $T^2\cdot x$. As in the proof of Lemma 6.4, $\Gamma$
acts on $T^2\cdot x$ with quotient a Klein bottle, since $\rho
(\Gamma)$ is not trivial. Therefore, the same argument in the
proof of Lemma 6.4 implies an element $t_0\in T^2$ so that
$t_0\gamma$ has a $3$-dimensional fixed point set $F$ in $\tilde
M$. By the Frankel's theorem, $F$ intersects with the two
$3$-dimensional strata, of circle isotropy groups. Clearly, $F$
projects to the fixed point set of $\Gamma$ in $D^3$. Recall that
$F$ intersects with a principal orbit in two circles. This
together shows that the fixed point set of $\Gamma$ in $D^3$ is a
$2$-dimensional, and so a disk, with non-empty intersections with
both $D^2_+$ and $D^2_-$. Therefore, the fixed point set
$(D^3)^\Gamma$ contains at least a point of $D^2_+\cap D^2_-$,
the fixed point set of the $T^2$-action. For any such a point
$[x]$, its preimage $x\in \tilde M$ satisfies $\gamma x=x$. A
contradiction, since $\Gamma$ acts freely on $\tilde M$.

Case 2. The $T^2$-action has no fixed point, but it has a
$3$-dimensional stratum with circle isotropy group.

Let $T^1\subset T^2$ denote the unique circle isotropy group with
$3$-dimensional fixed point set. Since $\Gamma$ preserves the
strata, $\Gamma$ preserves the isotropy group $T^1$, that is, for
any $g\in \Gamma$, $\rho (g)(T^1)=T^1$. Therefore, $T^1\rtimes
\Gamma$ acts on $\tilde M$. We claim that the {\it $T^1$-action
and $\Gamma$-action commute.} Then $M$ admits a $T^1$-action with
three dimensional fixed point set. By [GS] again we know that $M$
is diffeomorphic to a lens space.

It is clear that the $T^1$-action on $\tilde M$ is semi-free with
fixed point set a totally geodesic $S^3$, and the orbit space
$\tilde M^*$ is homeomorphic to $D^4$. Note that $\Gamma$ acts
freely on $\partial (\tilde M^*)=S^3$. Therefore $\Gamma$ acts
on $D^4$ with a unique fixed point in the interior, saying $0\in
D^4$, which is a principal orbit for the $T^1$-action on $\tilde
M$. If the commutativity fails, then there is a generator $\gamma
\in \Gamma$ such that $\rho (\gamma) \in \text{Aut} (T^1)$ is
given by the inverse automorphism. Let $T^1\cdot x_0$ denote the
principal $T^1$-orbit over $0\in D^4$. Since $\gamma (x_0)\ne
x_0$, let $t\in T^1$ satisfy the equation $t^2x_0=\gamma x_0$.
Then $x=tx_0$ satisfies the equation $\gamma x= \gamma tx_0=\rho
(\gamma)(t)\gamma x_0=t^{-1}\gamma x_0 =x$. A contradiction, since
$\gamma$ acts freely on $T^1\cdot x_0$.

Case 3. The $T^2$-action has only isolated singular orbits but
has finite order isotropy groups.

Assume that $\Bbb Z_p\subset T^2$ is an isotropy group of order
$p$, whose fixed point set $\tilde M^{\Bbb Z_p}$ is a totally
geodesic $3$-sphere. Note that $\Gamma$ preserves $\Bbb Z_p$ and
it acts freely on $\tilde M^{\Bbb Z_p}$. By Theorem 4.4 the
$T^2$-action has three isolated singular orbits, and the orbit
space $\tilde M^*$ is a homotopy $3$-sphere. We first claim that
{\it the $\Gamma$-action and the $T^2$-action commute} and thus
$T^2$ acts on $M$. Consequently, $\Gamma$ is cyclic. We argue by
contradiction. Suppose not, there is a nontrivial $\Gamma$-action
on the orbit space $\tilde M^*$ acting transitively on the set
of three singular orbits (by Lemma 6.5). Therefore, there is an
element $\gamma \in \Gamma$ which moves the three points
transitively. Note that the image of $\tilde M^{\Bbb Z_p}/T^2$ is
an interval, $[0,1]$, connecting two singular orbits. Let $\rho
:\Gamma \to \text{Aut}(T^2)$ denote the holonomy. By Lemma 6.5
$\rho (\gamma)$ is non-trivial. Note that $\Bbb Z_p$,
$\rho(\gamma)(\Bbb Z_p)$, $\rho (\gamma ^2)(\Bbb Z_p)$ are all
isotropy groups of the $T^2$-action on $\tilde M$ such that they
have pairwisely different fixed point set. However, since any two
of $\Bbb Z_p$, $\rho(\gamma )(\Bbb Z_p)$, $\rho (\gamma ^2)(\Bbb
Z_p)$ generate the same subgroup isomorphic to $\Bbb Z_p\oplus
\Bbb Z_p$, of rank $2$ in $T^2$, so the fixed point set of $\Bbb
Z_p\oplus \Bbb Z_p$ is the union of three isolated singular circle
orbits in $\tilde M\approx S^5$. A contradiction, by Theorem 4.3.

Consider the $T^2$-action on $M$. For a finite isotropy group
$\Bbb Z_p\subset T^2$, the $3$-dimensional fixed point component
$F\subset M^{\Bbb Z_p}$ is a lens space since it is invariant by
the isometric $T^2$-action (by [GS]). By Theorem 4.10 the
fundamental group of $F$ is isomorphic to $\Gamma$. Since the
orbit space $M/T^2$ is again a homotopy $3$-sphere with exactly
three isolated singular circle orbits, $M/T^2 -F/T^2$
homotopically retracts to the unique singular orbit (a point in
the orbit space) outside $F/T^2$. Thus, $M-F$ is homotopy
equivalent to the circle orbit outside $F$. This implies that
$M=D(\nu )\cup _\partial S^1\times D^4$, where $D(\nu)$ is the
normal disk bundle of $F$ in $M$, and $S^1\times D^4$ is a disk
tubular neighborhood of the circle orbit. By the Gysin exact
sequence for the unit circle bundle $S(\nu)$ over the lens space
$F=S^3/\Gamma$ we conclude that the Euler class $e(\nu )\in
H^2(F)\cong \Gamma $ is a generator. Thus, the bundle $D(\nu)$ is
uniquely up to conjugation.  Using the same fact about the gluing
along $S^1\times S^3$ as in the proof of Theorem 7.1 we conclude
that $M$ is diffeomorphic to a lens space, the desired result
follows. \qed\enddemo

\vskip4mm

\head 9. Proof of Theorem F
\endhead

\vskip4mm

 As we noticed in Section 3, Theorem F' implies Theorem F. The
 goal of this section is to prove Theorem F'. Recall that
the algebraic conditions on the fundamental group $\pi _1(M)$ in
Theorem F' are essentially as follows:

\noindent (9.1.1) Any index $\le 2$ normal subgroup $\Gamma
\triangleleft \pi_1(M)$ has a center $C(\Gamma)$ of index at least $ w$.

\noindent (9.1.2) Any index $\le 2$ normal subgroup $\Gamma
\triangleleft \pi_1(M)$ is not a spherical $3$-space group.

In fact (9.1.2) may be replaced by $\Gamma$ is neither cyclic, nor
generalized quaternionic group, and binary dihedral group.

\vskip 2mm

\proclaim{Lemma 9.2}

Let $M$ be a closed $5$-manifold of positive
sectional curvature which admits a $\pi_1$-invariant fixed point
free isometric $T^1$-action. If $\pi _1(M)$ satisfies (9.1.1) and
(9.1.2), then

\noindent (9.2.1) Every non-principal $T^1$-orbit is isolated, and
$\tilde M^*$ is a simply connected orbifold with only
isolated singularities.

\noindent (9.2.2) $H_2(\tilde M;\Bbb Z)$ is torsion free and has
rank equal to $b_2(\tilde M^*)-1$.
\endproclaim

\vskip2mm

\demo{Proof} (9.2.1) If there is a nontrivial
subgroup $\Bbb Z_p\subset T^1$ with a fixed point component
$\tilde M^{\Bbb Z_p}$ of dimension $3$, note that the
fixed point set $\tilde M^{\Bbb Z_p}$ is invariant by
the free $\pi _1(M)$-action. Thus $\pi_1(M)$ is isomorphic to the
fundamental group of a $3$-dimensional spherical space form, by
Theorem 4.10. A contradiction to the algebraic condition (9.1.2).

(9.2.2) First, $H_2(\tilde M;\Bbb Z)\cong \pi_2(\tilde M)$ (by the
Hurewicz theorem). Let $\tilde M_0$ denote the union of all
principal $T^1$-orbits. By (9.2.1), $\pi_i(\tilde M)\cong
\pi_i(\tilde M_0)$, $i=1, 2$. Because $\tilde M_0$ is obtained by
removing some isolated circle orbits, by the transversality
$\tilde M_0$ is simply connected and thus $\tilde M^*_0=\tilde M_0/T^2$
is simply connected. Since every singularity in $\tilde M^*$ is
a conical
point whose neighborhood in $\tilde M^*$ is a cone over a lens space
$S^3/\Bbb Z_{p}$, it is easy to see that $H_2(\tilde M^*;\Bbb
Z)=H_2(\tilde M^*_0;\Bbb Z)\cong \pi_2(\tilde M^*_0)$ (the last
isomorphism is from Hurewicz theorem). From the homotopy exact
sequence of the fibration,
$$1\to \pi_2(\tilde M_0)\to \pi_2(\tilde M^*_0)\to \Bbb Z\to 1,$$
it suffices to show that $H_2(\tilde M^*_0;\Bbb Z)$
is torsion free. This is true because $H_2(\tilde M^*_0;\Bbb Z)
\cong H^2(\tilde M^*_0,\partial
\tilde M^*_0;\Bbb Z)$ (the Lefschetz duality) whose torsion is
isomorphic to the torsion of $H_1(\tilde M^*_0,\partial
\tilde M^*_0;\Bbb Z)=0$ (the universal coefficient theorem).
\qed\enddemo

\vskip4mm

\subhead b. Estimate the Euler characteristic of $\tilde M^*$
\endsubhead

\vskip4mm

To determine the topology of $\tilde M$, we will first estimate
the Euler characteristic of $\tilde M^*$ using the constraint on
the fundamental groups.

\vskip2mm

\proclaim{Proposition 9.3}

Let the assumptions and notations be as
in Lemma 9.2. Then $\chi(\tilde M^*) =2+b_2(\tilde M^*)\le 5$ or
equivalently, $b_2(\tilde M^*)=b_2(\tilde M)+1\le 3$.
\endproclaim

\vskip2mm

Let $p: \tilde M\to \tilde M^*$ denote the orbit
projection. Then $\pi_1(M)$ acts on $\tilde M^*$ by isometries. For
$\gamma\in \pi_1(M)$, we will use $\gamma^*$ to denote the
isometry on $\tilde M^*$ induced by $\gamma$ (cf. $\S$2). By Theorem
4.7, $(\tilde M^*)^{\gamma^*}\ne \emptyset$, if $\gamma$ commutes
with the $T^1$-action. To estimate the characteristic of $\tilde
M^*$, we will first estimate the number of isolated $
\gamma^*$-fixed points (see Theorem 4.1) via the technique of
$q$-extent estimate ([GM], [Ya]).

The $q$-extent $xt_q(X)$, $q\ge 2$, of a compact metric space $(X, d)$ is, by
definition, given by the following formula:
$$
xt_q(X)={q \choose 2} ^{-1}\text{max} \Bigl\{ \sum _{1\le i<j\le q}
d(x_i, x_j): \{x_i\}_{i=1}^q \subset X \Bigr\}
$$

Given a positive integer $n$ and integers $k, l \in \Bbb Z$ coprime to $n$,
let $L(n; k, l)$ be the $3$-dimensional lens space, the quotient space of a
free isometric $\Bbb Z_n$-action on $S^3$  defined by
$$\psi _{k,l}: \Bbb Z_n \times S^3\to S^3; \text{  } g(z_1, z_2)=
(\omega ^kz_1,\omega ^l z_2 )$$
with $g\in \Bbb Z_n$ a generator, $\omega =e^{i\frac{2\pi}n}$ and
$(z_1, z_2) \in S^3\subset \Bbb C^2$.

Note that $L(n;k, l)$ and $L(n; -k, l)$ (resp. $L(n; l, k)$) are
isometric (cf. [Ya] p.536). Obviously $L(n; -k,l)$ and $L(n;
n-k,l)$ are isometric. Therefore, up to isometry we may {\it
always} assume $k, l\in (0, n/2)$ without loss of generality. The
proof of Lemma 7.3 in [Ya] works identically for $L(n; k, l)$ with
$0<k, l<n/2$ to prove

\vskip2mm

\proclaim{Lemma 9.4 ([Ya])}

Let $L(n; k, l)$ be a $3$-dimensional
lens space of constant sectional curvature one. Then
$$\aligned
xt_q(L(n;k,l))&\le \text{arccos}\Bigl\{ \text{cos}
(\alpha _q)\text{cos }\pi n^{-\frac 12}\\
&-\frac 12 \{ (\text{cos }\pi n^{-\frac 12}-\text{cos }\pi /n)^2+
\text{sin}^2(\alpha _q) (n^\frac 12 \text{sin }\pi/n -\text{sin }
\pi  n^{-\frac 12})^2 \}^ {\frac 12} \Bigr\}\endaligned
$$
where $\alpha _q=\pi /(2(2-[(q+1)/2]^{-1}))$, where ``$[\,x\,]$''
means the integer part of $x$.
\endproclaim

\vskip2mm

\proclaim{Corollary 9.5}

Let $L(n;k,l)$ be a $3$-dimensional lens
space of constant sectional curvature one. If $n\ge 61$, then
$xt_5(L(n; k, l))< \pi/3$.
\endproclaim

\vskip2mm

\proclaim{Corollary 9.6}

If the exponent of $\gamma^*$ is at
least $61$, then $(\tilde M^*)^{\gamma^*}$ contains at most five
isolated fixed points.
\endproclaim

\vskip2mm

\demo{Proof} We argue by contradiction, assuming $x^*_1,...,
x^*_6$ are six isolated $\gamma^*$-fixed points. Let
$X=\tilde M^*/\langle \gamma^*\rangle$. Connecting each pair of
points by a minimal geodesic in $X$, we obtain a
configuration consisting of twenty geodesic triangles. Because
$X$ has positive curvature in the comparison sense ([Pe]),
the sum of the interior angles of each triangle is $>\pi$ and thus
the sum of total angles of the twenty triangles, $\sum
\theta_i>20\pi$. We then estimate the sum of the total angles in
the following way, first estimate from above of the ten angles
around each $x^*_i$ and then sum up over the six points. We
claim that the sum of angles at $x^*_i$ is bounded above by
$10\cdot xt_5(x^*_i)\le 10\frac {\pi}3$ and thus $\sum
\theta_i\le 6(10\cdot \frac{\pi}3)=20\pi$, a contradiction.

Let $\tilde x_i\in \tilde M$ such that $p(\tilde x_i)=x^*_i$,
and let $\tilde t\in T^1$ such that $\tilde t\cdot \gamma$ fixes
$T^1 (\tilde x_i)$. Let $S^\perp_{\tilde x_i}$ denote the unit
$3$-sphere in the normal space of $T^1(\tilde x_i)$. If the
isotropy group at $\tilde x_i$ is trivial, then the space of
directions at $x^*_i$ in $X$ is isometric $S^\perp_{\tilde
x_i}/\langle \tilde t\cdot \gamma\rangle$ which is a lens space
with a fundamental group of order $|\gamma^*|$. By Corollary
9.4, we conclude that the sum of the ten angles is bounded above
by
$$\left(\matrix 5\\2\endmatrix\right)xt_5(L)=10\cdot \frac{\pi}3.$$
If the isotropy group at $\tilde x_i$ is not trivial, then the
above estimate still holds because the $5$-extent only gets
smaller when passing to the quotient of $S^\perp_{\tilde x_i}$ by
the isotropy group. \qed\enddemo

\vskip2mm

\proclaim{Lemma 9.7}

Let $\tilde X$ be a simply connected
topological space of dimension $4$ with a finite group $G$-action.
Assume that the total Betti number of $\tilde X$ is bounded above by a
constant $N$. Then $G$ has a normal subgroup of order at least
$|G|/c(N)$ which acts trivially on the homology $H^*(\tilde X)$,
where $c(N)$ is a function depending only on $N$.
\endproclaim

\vskip2mm

\demo{Proof} Let $\rho: G\to \text{Aut}(H^*(\tilde X))$ denote the
homomorphism induced by the $G$-action on $H^*(\tilde X)$. Because
the torsion subgroup of $\text{Aut}(H^*(\tilde X))$ is a subgroup
of $\text{GL}(\Bbb Z, N)$, which is bounded above by a constant
depending only on $N$ (cf. [Th]), say $c(N)$, the conclusion
follows. \qed\enddemo

\vskip2mm

\demo{Proof of Proposition 9.3}

Because $\tilde M^*$ is a simply connected orbifold with isolated
singularities (Lemma 9.2), we can apply the Poincar\'e duality on
homology groups with rational coefficients to conclude that
$\chi(\tilde M^*)=2+b_2(\tilde M^*)$. Let $\Gamma _0$ be the principal
isotropy group of $\pi _1(M)$ on $\tilde M^*$. By Theorem 4.9
$\Gamma _0$ is cyclic and belongs to the center of a certain index
at most two normal subgroup (the subgroup of orientation
preserving isometries in $\pi _1(M)$). By the assumption (9.1.1),
the index $[\pi _1(M):\Gamma _0]\ge w$.

On the other hand, by the Betti number bound in [Gr] and Lemma 9.7
we may assume that the normal subgroup $G \subset \pi _1(M)/\Gamma _0$
has order $\ge w/k_0$, and thus acts trivially on
$H^*(\tilde M^*;\Bbb Z)$.

For any $\beta^* \in G$, by Theorem 4.8 the fixed point set
$(\tilde M^*)^{\beta^*}\ne \emptyset$ and by Theorem 4.10 we may
assume that $(\tilde M^*)^{\beta^*}$ is a finite set (otherwise
$\tilde M$ contains a $3$-dimensional totally geodesic
submanifold, and so Theorem 4.10 implies that it is a homotopy
sphere). If $\beta^*$ is of prime order and $|\beta^* |\ge
61$, by Theorem 4.1
$$\chi(\tilde M^*)=\chi(F(\beta^* ,\tilde M^*)).$$
Therefore, by Corollary 9.6 we conclude $\chi (\tilde M^*)\le 5$.

Now we assume that $|G|$ has all prime factors $\le 60$. Since
there are at most $17$ primes less than $60$, there is a prime
$p\le 60$ so that $G$ has a $p$-Sylow subgroup $G_p$ of order $\ge
\frac{w}{17k_0}$. Let $\beta^*_0\in G_p$ be an element of order
$p$ generating a normal subgroup of $G_p$ (by finite group theory
[Se], $G_p$ is nilpotent), we may assume that $\beta^* _0$ has
only isolated fixed points for the same reasoning as above, say
$p_1, \cdots, p_n$. By Theorem 4.1, $n=\chi ((\tilde M^*)^{\beta^*_0})=\chi (\tilde M^*)$. Observe that $n\le b=b(4)$. Now $G_p$
acts on the fixed point set $(\tilde M^*)^{\beta^*_0}$, thus we
get a homomorphism $h: G_p\to S_n$, where $S_n$ is the permutation
group of $n$-words. Therefore, the kernel of $h$ has order at
least $\frac w{17k_0\cdot n!}\ge \frac w{17k_0\cdot b!}$. By a
well-known result of Gromov $\pi _1(M)$ may be generated by a
bounded number of generators, so is $G_p$, generated by $c$
elements. Hence $\text{ker}(h)$ contains an element, say $
\beta^* _1$, of order at least $61$, if $\frac w{17k_0\cdot b!}$ is
sufficiently large. Since $\beta^* _1$ fixes the set $\{p_1,
\cdots, p_n\}$ pointwisely, by Corollary 9.6 $n\le 5$. Therefore,
$\chi (\tilde M^*)=n\le 5$. \qed\enddemo

\vskip4mm

\subhead c. The completion of the proof of Theorem F'
\endsubhead

\vskip4mm

\proclaim{Lemma 9.8}

Let $M$ be a closed $5$-manifold of positive
sectional curvature which admits a $\pi _1$-invariant isometric
$T^1$-action. If $\pi_1(M)$ satisfies (9.1.1), then the
$T^1$-action on $\tilde M$ is free.
\endproclaim

\vskip2mm

\demo{Proof} We argue by contradiction, assuming the $T^1$-action
is not free. Then there is at least a finite isotropy group $\Bbb
Z_p\subset T^1$. By Theorem 4.10 we may assume that $\dim(F(\Bbb
Z_p,\tilde M))=1$, since otherwise, $\tilde M$ contains a totally
geodesic $3$-manifold. Then $\tilde M^{\Bbb Z_p}$ consists of at
most two components (circles), if $b_2(\tilde M)\le 1$, or three
components if $b_2(\tilde M)\le 2$ (Theorem 4.2).

Let $H$ denote the subgroup of $\pi_1(M)$ preserving all
components of $\tilde M^{\Bbb Z_p}$. Then $H$ is a cyclic normal
subgroup such that the quotient $\pi _1(M)/H$ acting effectively
on the set of exceptional orbits. If $b_2(\tilde M)\le 1$, by
counting the number of components $\pi _1(M)/H$ has order at most
$2$.  A contradiction to (9.1.1).

If $b_2(\tilde M)=2$, and $\tilde M^{\Bbb Z_p}$ has exactly three
components, it is easy to see that $H$ acts on $\tilde M^*$ has at most five isolated fixed points, three of them are
the exceptional $T^1$-orbits with isotropy group $\Bbb Z_p$.
Because $H$ is normal in $\pi _1(M)$, thus $\pi _1(M)$ acts on the
fixed point set and it sends an exceptional orbit to an
exceptional orbit. Therefore, $\pi _1(M)$ acts on the union of the
rest at most two isolated $T^1$-orbits fixed by $H$. This implies
once again that $\pi_1(M)$ has an index $\le 2$ cyclic normal
subgroup. The desired result follows.
\qed\enddemo

\vskip2mm

\proclaim{Lemma 9.9}

Let $M$ be as in Lemma 9.8. Then the induced
$\pi_1(M)$-action on $\tilde M^*$ is pseudo-free.
\endproclaim

\vskip2mm

\demo{Proof} By Lemma 9.8, $\tilde M^*$ is a closed
$4$-manifold. If there is an element $\gamma\in \pi_1(M)$ with a
fixed point component $F\subset \tilde M^*$ of dimension $2$, there
exists an element $t\in T^1$ so that $t\gamma$ has a
$3$-dimensional fixed point set in $\tilde M$. By Theorem 4.10
this implies that $\tilde M\approx S^5$. By Lemma 1.10, we conclude
that $\pi _1(M)$ acts on $\tilde M^*$ pseudofreely, a contradiction.
\qed\enddemo

\vskip2mm

\proclaim{Lemma 9.10}

Let $M$ be a closed $5$-manifold of positive
sectional curvature which admits a $\pi_1$-invariant fixed point
free isometric $T^1$-action. If $\pi _1(M)$ satisfies (9.1.1) and
(9.1.2), then $\tilde M$ is a homotopy sphere.
\endproclaim

\vskip2mm

\demo{Proof} Because the $T^1$-action on $\tilde M$ is free (Lemma
9.7), $\tilde M^*$ is a closed simply connected smooth $4$-manifold
of positive sectional curvature. Because $b_2(\tilde M^*)=b_2
(\tilde M)+1\ge 1$, by Proposition 9.3, $3\le 2+b_2(\tilde
M^*)=\chi(\tilde M^*) \le 5$.

If $\chi(\tilde M^*)=3$, then $b_2(\tilde M)=b_2(\tilde M^*)-1=0$.
This together with (9.2.2) implies that  $\tilde M\approx S^5$.

It suffices to rule out the cases $\chi(\tilde M^*)=4$ and $\chi(\tilde M^*)=5$.

Case (a) If $\chi(\tilde M^*)=4$;

Note that $\tilde M^*$ is homeomorphic to
$S^2\times S^2$, $\Bbb CP^2\#\Bbb  CP^2$ or $\Bbb CP^2\#
\overline{\Bbb CP^2}$, up to a possible orientation reversing
([Fr]). By the classification of simply connected $5$-manifolds
(cf. [Ba]) $\tilde M=S^2\times S^3$ or $S^2\tilde \times S^3$.
Because $\pi_1(M)$ preserves the Euler class of the principal
circle bundle $T^1\to \tilde M\to \tilde M^*$, the natural
homomorphism given by the $\pi_1(M)$-action on the second homology
group, $\alpha: \pi _1(M)\to \text{Aut}(H^2(\tilde M^*;\Bbb Z);
I)$, where $\text{Aut}(H^2(\tilde M^*;\Bbb Z); I)$ is the
automorphism group preserving the intersection form.

Let $\Gamma _0$ (resp. $\Gamma$) be the principal isotropy group
of $\pi _1(M)$-action on $\tilde M^*$ (resp. orientation preserving
subgroup of $\pi _1(M)$). Recall that $\Gamma\subset \pi _1(M)$ is
a normal subgroup of index at most $2$, and $\Gamma _0$ is in the
center of $\Gamma$ which generates a cyclic subgroup with any
element of $\Gamma$ (cf. 4.10) Let $G=\Gamma /\Gamma _0$. The
homomorphism $\alpha$ reduces to a homomorphism $\bar \alpha :
G\to \text{Aut}(H^2(\tilde M^*;\Bbb Z); I)$. It is easy to see that
the image of $\bar \alpha$ has order at most $2$ (cf. [Mc2]).

Subcase (a1) If $\text{ker}(\bar \alpha)$ has odd order;

Consider the action of $\text{ker}(\bar \alpha)$ on $\tilde M^*$.
Observe that the action is pseudo-free, and every isotropy group
is cyclic. By [Mc2] Lemma 7.5 and the first paragraph in the proof
of Lemma 4.7 [Mc2] (which identically extends to our case) it
follows that, $\text{ker}(\alpha^*)$ is cyclic of odd order. If
$Im (\alpha^*)$ is not trivial and $\tilde M^*\approx S^2\times
S^2$, by [Bre] VII Corollary 7.5 there is an involution on $\tilde M^*$
with a $2$-dimensional fixed point set. A contradiction by Lemma
9.8. If $\tilde M^*\approx \Bbb CP^2\# \Bbb CP^2$ or $\Bbb CP^2\#
\overline{\Bbb CP}^2$, in the former case every self homeomorphism
of $\tilde M^*$ is orientation preserving, and for the latter
$Im(\alpha^*)=0$. Therefore, in either cases, $\pi _1(M)$ has
an image in $\text{Aut}(H^2 (\tilde M^*;\Bbb Z); I)$ of order at
most $2$, and by Corollary 4.10 it contains a normal cyclic
subgroup of index at most $2$, a contradiction to (9.1.1).

Subcase (a2) If $\text{ker}(\alpha^*)$ has even order;

By [Bre] VII Lemma 7.4, any involution in $\text{ker}(\alpha^*)$
has a $2$-dimensional fixed point set on $\Bbb CP^2\#\Bbb CP^2$ or
$\Bbb CP^2\#\overline{ \Bbb CP}^2$. Therefore, by Lemma 9.9 we may
assume that $\tilde M^*=S^2\times S^2$. By [Mc2] Theorem 3.3,
$\text{ker}(\alpha^*)$  is a polyhedral. Therefore
$\text{ker}(\alpha^*)$ is either cyclic, dihedral, or is a
non-abelian group of order of order $12$ (Tetrahedral
group and two others) or
$\text{ker}(\alpha^*)$ is Octahedral group (of order $24$), or
Icosahedral group (of order $60$). Thus, $\text{ker}(\alpha^*)$ is cyclic or dihedral
if the order $|G|>120$. We may require our constant $w$ in Theorem
E is larger than $120$. By [Mc2] Theorems 3.9 and 4.10 $G$ is
listed as follows:

\noindent (9.10.1) a cyclic, or dihedral group;

\noindent (9.10.2) $Q_{2^km}\times \Bbb Z_n$; where $Q_{2^km}$ is the
generalized quaternionic group, $m, n$ are coprime odd integers;

\noindent (9.10.3) $D_{2^km}\times \Bbb Z_n$, where $D_{2^km}=\Bbb Z_{m}
\rtimes \Bbb Z_{2^k}$ and $\Bbb Z_{2^k}$ acts on $\Bbb Z_m$ by
inverse automorphism, $m, n$ are coprime odd integers, and $k\ge
2$;

\noindent (9.10.4) A non-splitting extension of a dihedral group
by $\Bbb Z_2$.

Since $\Gamma$ is a center extension of a cyclic group by $G$, by
Lemma 9.11 below the group $\Gamma$ satisfies $2p$-condition,
i.e., for any prime $p$, a subgroup of order $2p$ is cyclic (cf.
[Mi]). By group extension theory, it is not hard to verify, the
dihedral groups in (9.10.1), (9.10.2) and (9.10.4) must be reduced
from a quaternionic subgroup of $\Gamma$, and the group $D_{2^km}$
in (9.10.3) is reduced from $D_{2^{k'}m}$, where $k'\ge 2$.
Moreover, for $G$ in (9.10.1), (9.10.2) and (9.10.3), $\Gamma$ is
isomorphic to a $3$-dimensional spherical space form group
(compare [Mi] Theorem 2), and so $\pi_1(M)$ contains an index at
most $2$ subgroup isomorphic to a $3$-dimensional spherical space
form group, a contradiction by the assumption. If $G$ is a group
in (9.10.4), and $[\pi _1(M): \Gamma]=2$, there is a {\it maximal}
cyclic subgroup $\langle \gamma \rangle$ of $\pi _1(M)$ of index
$8$. Therefore, $\pi _1(M)$ acts on the $4$ fixed points of
$\gamma$ on $S^2\times S^2$ without any even order isotropy group
by the maximality of $\langle \gamma \rangle$. A contradiction,
since the permutation group $S_4$ does not contain cyclic subgroup
of order $8$. This shows that $\Gamma =\pi _1(M)$, and it contains
a quaternionic subgroup of index $2$. A contradiction again to the
assumption.

Case (b) If $\chi (\tilde M^*)=5$;

Since the Euler characteristic $\chi (\tilde M^*)$ is odd, by [Bre]
VII Corollary 7.6 any involution on $X$ has a $2$-dimensional
fixed point set. Therefore, by Lemma 9.9 we may assume that $\pi
_1(M)/\Gamma _0$ has odd order. In particular, $\pi _1(M)=\Gamma$,
i.e. any element of $\pi _1(M)$ preserves the orientation of
$\tilde M^*$. For the homomorphism $\alpha^* : G\to
\text{Aut}(H_2 (\tilde M^*;\Bbb Z))$, by [Mc1] the kernel
$G_0=\text{ker}(\bar \alpha)$ is an abelian subgroup of $T^2$ with
non-empty fixed point set. Therefore, $G_0$ must be cyclic by
combining Lemma 9.9, otherwise, there exists a cyclic subgroup of
$G_0$ with a $2$-dimensional fixed point set. Therefore $\pi
_1(M)/G_0$ is isomorphic to an odd order subgroup of the
automorphism group $\text{Aut}(H_2(\tilde M^*;\Bbb Z))=
\text{GL}(\Bbb Z, 3)$. By [Th] the finite subgroup of
$\text{GL}(\Bbb Z, 3)$, up to possible $2$-torsion, is a subgroup
of $\text{GL} (\Bbb Z_2, 3)$ which has order
$(2^3-1)(2^3-2)(2^3-4)$, which is coprime to $5$. Note that $\pi
_1(M)$ preserves the fixed point set $(\tilde M^*)^{G_0}$,
which consists of $5$ isolated points. Therefore $\pi _1(M)/G_0$
is an odd order subgroup of $S_5$ and $\pi _1(M)/ G_0\ne \Bbb
Z_5$. This implies that $\pi _1(M)$ fixes at least $2$ points
among $(\tilde M^*)^{G_0}$ and so $\pi _1(M)$ is cyclic.  A
contradiction.
\qed\enddemo

\vskip2mm

\proclaim{Lemma 9.11}

Let $S^2\times S^3\to S^2\times S^2$ be a $G$-equivariant
principal $T^1$-bundle. If $G$ acts freely on $S^2\times S^3$,
pseudofreely on $S^2\times S^2$. Then there is no order $2$
element of $G$ acting non-trivially on $S^2\times S^2$ inducing
the identity on homology groups.
\endproclaim

\vskip2mm

\demo{Proof} If not, for an order $2$ element $\beta$, its fixed
point set $F$ in $S^2\times S^2$ consists of $4$ points (by
Theorem 4.1).
 For such a fixed point $[x]\in S^2\times S^2$ with
$x\in S^2\times S^3$, there is an element $t\in T^1$ so that
$\beta x=tx$. The freeness of $T^1$-action and $\beta$-action
implies that the order of $t$ is the same as $\beta$, of order
$2$. Since $T^1$ has only one element of order $2$, say $t_0$,
this proves that $\beta t_0^{-1}$ has fixed point set the union of
$4$ circles. A contradiction by Theorem 4.2. \qed\enddemo

\vskip14mm

\head Appendix: $s$-cobordism theory
\endhead

\vskip4mm

In the proof of Theorem E for pseudofree $T^2$-action (Section 7)
we assumed the Poincar\'e conjecture holds. The goal of this
appendix is to give an account how this assumption can be removed
by using the $s$-cobordism theory.

We call two closed $n$-manifolds $M_1$ and $M_2$ are {\it
$s$-cobordant} if there is an $(n+1)$-manifold $W$ with boundary
$M_1\sqcup M_2$ so that the inclusions $i_1: M_1\to W$ and $i_2:
M_2\to W$ are both simple homotopy equivalences (cf. [Ke]). The
manifold $W$ is called an {\it $s$-cobordism}. The deep
$s$-cobordism theorem (originally due to Smale for simply
connected manifold, and extended to non-simply connected case by
Barden-Mazur-Stallings, cf. [Ke] for a detailed account) asserts
that an $s$-cobordism is diffeomorphic (or homeomorphic if
manifolds are TOP.) to the product $M_1\times [0,1]$, provided
$n\ge 5$. Hence, $M_1$ and $M_2$ are diffeomorphic (resp.
homeomorphic). The dimension assumption is crucial. In fact,
counterexamples exist for $n=3, 4$, by Cappell-Shaneson, and
Donaldson's theory. An $s$-cobordism theorem for simply connected
$4$-manifold (i.e. $n=3$) would imply the $3$-dimensional
Poincar\'e conjecture.

It is a well-known fact that every simply connected $3$-manifold
is $s$-cobordant (by terminology should be called $h$-cobordant)
to $S^3$. The key issue of our solution in the proof Theorem E
without assuming the Poincar\'e conjecture is to replace the
homeomorphism (equivalently diffeomorphism) by $s$-cobordism. If
we could obtain an $s$-cobordism between our $5$-manifold $M$ with
a spherical space form, the $s$-cobordism theorem applies to imply

our desired result. Note that the dimension shifting is important
and this can be obtained by using  the additional $\pi
_1$-invariant $T^2$-action on the manifold.

\vskip2mm

\definition{Definition A.1}

Let $M_1$ (resp. $M_2$) be a closed $5$-manifold with a smooth
action by $T^2\rtimes G$ where the $T^2$-action is pseudofree and
the $G$-action is free. We call that $W$ is an $s$-cobordism
between $(M_1, T^2\rtimes G)$ and $(M_2, T^2\rtimes G)$ if

\noindent (A.1.1) $W$ is an $s$-cobordism between $M_1$ and $M_2$;

\noindent (A.1.2) There is an action by $T^2\rtimes G$ on $W$ whose
restriction on $M_1$ (resp. $M_2$) gives $(M_1, T^2\rtimes
G)$(resp. $(M_1, T^2\rtimes G)$ ), and the $T^2$-action is
pseudofree and the $G$-action is free.
\enddefinition

\vskip2mm

By the $s$-cobordism theory (cf. [Ke]), $W$ is homeomorphic to the
product $M_i\times [0,1]$ when $\dim(M_i)\ge 5$. However, the
induced $T^2$-action on $M_i\times [0,1]$ may not preserve the
slice $M_i\times \{t\}$ ($i=1, 2$), $t\in (0,1)$, so that one may
not conclude that $(M_1,T^2)$ is conjugate to $(M_2,T^2)$.

\vskip2mm

As we have seen before, if $T^2$ acts pseudofree on $S^5$, the
orbit space $S^5/T^2$ is a manifold homotopy equivalent to $S^3$.
Hence, if $W$ is an $s$-cobordism between two pseudofree
$T^2$-actions on $S^5$, obviously, the orbit space $W/T^2$ gives
an $s$-cobordism between the two orbit spaces (two homotopy
$3$-spheres) of the actions on $S^5$. We will show below the
converse also holds.

\vskip2mm

\proclaim{Lemma A.2}

A pseudo-free $T^2$-action on $S^5$ is $s$-cobordant to a
pseudo-free linear $T^2$-action.
\endproclaim

\vskip2mm

\demo{Proof} As in Section 7 we make use the convention that
$T^2_0$ indicates a linear $T^2$-action on $S^5$. Let $V$ be an
$s$-cobordism of the homotopy $3$-sphere $S^5/T^2:=\Sigma$ and
$S^5/T^2_0=S^3$. By Theorem 4.4 $T^2$ (resp. $T^2_0$) has three
isolated circle orbits, denoted by $T^1_i\times D^4$ (resp. $\hat
T^1_i\times D^4$) respectively. Let $p_i$ (resp. $q_i$) denote the
projections of these circle orbits in $S^5/T^2$ (resp.
$S^5/T^2_0$), $i=1, 2,3$, respectively. Connecting $p_i$ and $q_i$
by a simple disjoint paths in $V$, $i=1, 2, 3$, let $V_0$ denote
the complement of the three simple paths. Observe that $V_0$ is
homotopy equivalent to $S^3_0=S^3-q_1\cup q_2\cup q_3$ (resp.
$\tilde M^*_0=\tilde M^*_0-p_1\cup p_2\cup p_3$). It is easy to see that
$H^2(S^3_0;\Bbb Z^2)\cong \text{Hom}(\Bbb Z^2, \Bbb Z^2)$.
Consider the Euler classes of principal $T^2$-bundles on $\Sigma
_0$ and $S^3_0$ as $2\times 2$ matrices in the above group. Hence,
the determinants of the matrices are $\pm 1$ because the
$2$-connectedness of the complement of the singular orbits.
Therefore, modifying by an automorphism of $T^2$ if necessary (in
other words, up to weak equivalence), we may assume that the two
Euler classes are homotopic, regarded as maps to the classifying
space $BT^2$. The homotopy gives a map $h: V_0\to BT^2$ whose
restrictions on $S^3_0$ and $\tilde M^*_0$ correspond to the
classifying maps of the principal $T^2$-bundles. Let $W_0$ denote
the total space of the $T^2$-bundle over $V_0$, which is a
$6$-manifold with boundary.

Finally, we can attach equivariantly three copies of $T^1\times
D^4\times [0,1]$ to $W_0$ along the boundaries, so that the three
copies of $T^1\times D^4\times \{0\}$ fill the three singular
orbits in $S^5$. This gives a $T^2$-manifold $W$ with boundary
$S^5\sqcup S^5$, which yields the desired $s$-cobordism between
$(S^5,T^2)$ and $(S^5,T^2_0)$. \qed\enddemo

\vskip2mm

In Lemmas 7.3 and 7.4 the induced $\Gamma$-action on the homotopy
$3$-sphere $\Sigma$ (the orbit space) may not be trivial. We need
to find an equivariant $s$-cobordism between $(\Sigma ,\Gamma)$
with a standard linear action on $(S^3, \Gamma)$. Recently, a deep
result (cf. [BLP]) on $3$-orbifold implies that every finite group
acting non-freely on $S^3$ is conjugate to a linear action. Let
$\Gamma$ be a finite group. For a smooth or locally linear nonfree

$\Gamma$-action on a homotopy $3$-sphere (possibly reducible)
$\Sigma$, the result of [BLP] implies that $(\Sigma, \Gamma)$ is
$\Gamma$-equivariantly diffeomorphic to $ (S^3, \Gamma)\# (|\Gamma
|\cdot \tilde M^*_0, \Gamma)$, where $\Gamma$ acts linearly on $S^3$,
and acts freely on $|\Gamma |\cdot \tilde M^*_0$ and $\tilde M^*_0$ is a
homotopy $3$-sphere (we thank Porti for pointing out this fact to
us). Therefore, it is easy to get a $4$-dimensional $s$-cobordism
$(V, \Gamma)$ between $(\Sigma, \Gamma)$ and $(S^3, \Gamma)$.
Indeed, we may take $(V, \Gamma)=(S^3\times [0,1], \Gamma)
\natural (|\Gamma |\cdot B_0, \Gamma)$, where $B_0$ is a
contractible $4$-manifold with boundary $\tilde M^*_0$, and
$\natural$ is the boundary connected sum along the boundary piece
$S^3\times \{1\}$. An exceptional case is $\Gamma =\Bbb Z_3$ which
acts freely on the homotopy $3$-sphere $\Sigma$. In this case the
main result of [BLP] does not apply. However, it is easy to see
that $\Sigma /\Bbb Z_3$ is simple homotopy equivalent to the lens
space $S^3/\Bbb Z_3$ (unique), because that the Whitehead torsion
$Wh(\Bbb Z_3)=0$ (cf. [Coh]). Hence, there is an $s$-cobordism
between $S^3/\Bbb Z_3$ and $\Sigma /\Bbb Z_3$, which gives exactly
a $\Bbb Z_3$-equivariant  $s$-cobordism between $(S^3, \Bbb Z_3)$
and $(\Sigma, \Bbb Z_3)$.

\vskip2mm

\proclaim{Lemma A.3}

The $\Gamma$-equivariant principal $T^2$-bundle over $\tilde M^*_0$
in Lemmas 7.3 and 7.4 is $\Gamma$-equivariantly $s$-cobordant to a
linear principal $T^2$-bundle over $S^3_0$. Hence, $M_0$ is
diffeomorphic to the linear model $S^5_0/\Gamma _0$.
\endproclaim

\vskip2mm

\demo{Proof} Because the proofs of Lemmas 7.3 and 7.4 involve only
homotopy theory, using the $\Gamma$-equivariant $s$-cobordism $V$
defined above, there is a well-defined $\Gamma$-equivariant
principal $T^2$-bundle over $V_0$, where $V_0$ is obtained by
removing certain $\Gamma$-equivariant three simple paths, e.g., in
the case $\Gamma$ acts non-freely on $\Sigma$, just take the
product $(p_1\cup p_2\cup p_3)\times [0,1] \subset V$; in the case
$\Gamma$ acts freely on $\Sigma$, take the preimage of any simple
path joining the singular point of $S^3/\Bbb Z_3$ and $\Sigma
/\Bbb Z_3$. The total space of the principal $T^2$-bundle over
$V_0$ gives a $\Gamma$-equivariant $s$-cobordism between $(\tilde
M_0, \Gamma)$ and $(S^5_0, \Gamma _0)$, which implies the
diffeomorphism between $S^5_0/\Gamma _0$ and $M_0$. \qed\enddemo

\enddocument